\documentclass[12pt,reqno]{amsart}
\usepackage{amssymb,delarray}
\usepackage{amsfonts}
\usepackage{epsfig}
\usepackage[all]{xy}

\makeindex{}

\def\ge{\geqslant}

\def\pr{{\sf pr}}

\def\pr{\mathsf{pr}}

\newtheorem{thm}{Theorem}[section]
\newtheorem{lem}[thm]{Lemma}
\newtheorem{prop}[thm]{Proposition}
\newtheorem{claim}[thm]{Claim}
\newtheorem{df}[thm]{Definition}
\newtheorem{cor}[thm]{Corollary}

{\catcode`\@=11
\gdef\n@te#1#2{\leavevmode\vadjust{%
 {\setbox\z@\hbox to\z@{\strut#1}%
  \setbox\z@\hbox{\raise\dp\strutbox\box\z@}\ht\z@=\z@\dp\z@=\z@%
  #2\box\z@}}}
\gdef\leftnote#1{\n@te{\hss#1\quad}{}}
\gdef\rightnote#1{\n@te{\quad\kern-\leftskip#1\hss}{\moveright\hsize}}
\gdef\?{\FN@\qumark}
\gdef\qumark{\ifx\next"\DN@"##1"{\leftnote{\rm##1}}\else
 \DN@{\leftnote{\rm??}}\fi{\rm??}\next@}}

\begin{document}
\baselineskip=14.pt plus 2pt %13.7pt plus 2pt

\title[Doubling formula]{Full twists factorization
formula for double number of strings}
\author[Vik.S.~Kulikov]{Vik. S.~Kulikov}
\address{Steklov Mathematical Institute\\
Gubkina str., 8\\
119991 Moscow \\
Russia} \email{kulikov@mi.ras.ru}

%\curraddr{}
\dedicatory{} \subjclass{}
\thanks{The work  was partially supported by
RFBR ({\rm No.} 02-01-00786) and {\rm INTAS (No.} 00-0269). }
\date{This version: July 2003. }
\begin{abstract}
A formula for factorizations of the full twist in the braid group
$Br_{2m}$ depending on any four factorizations of  the full twist
in $Br_{m}$ is given. Applying this formula, a symplectic
4-manifold $X$ and two isotopic
%regular homotopic
generic
coverings $f_i :X\to \mathbb C\mathbb P^2$, branched,
respectively, along cuspidal Hurwitz curves $\bar H_i\subset
\mathbb C\mathbb P^2$ (without negative nodes) having different
braid monodromy factorization types, are constructed. The class of
the fundamental groups of the complements of the affine plane
Hurwitz curves is described in terms of generators and defining
relations.
\end{abstract}

\maketitle
%\pagenumbering{roman}
\setcounter{tocdepth}{2}
%\setcounter{tocdepth}{1}
%\tableofcontentsSymplectic isotopy of sections

%\setcounter{section}{6}

\def\st{{\sf st}}

\setcounter{section}{-1}
\section{Introduction}

Let $(X,L)$ be a polarized projective surface over the field
$\mathbb C$, where the polarization $L$ is a very ample line
bundle on $X$. Three generic sections of $L$ define a generic
covering $f:X\to \mathbb C\mathbb P^2$ branched over an algebraic
cuspidal curve $\bar H\subset \mathbb C\mathbb P^2$ (possibly, one
should take three sections of $L^{\otimes 2}$ instead of the
sections of $L$ to obtain a generic covering of the plane
(\cite{Ku-Ku})). Recently, Auroux and Katzarkov (see \cite{Au},
\cite{Au-Ka}) obtained
%proved
a similar result in symplectic case. More precisely, let
$(X,\omega)$ be a compact symplectic 4-manifold with symplectic
form $\omega$ whose class $[\omega]\in H^2(X,\mathbb Z)$. Fix an
$\omega$-compatible almost complex structure $J$ and corresponding
Riemannian metric $g$. Let $L$ be a line bundle on $X$ whose first
Chern class is $[\omega]$. Then, for $k>>0$, the line bundle
$L^{\otimes k}$ admits many approximately holomorphic sections so
that one can choose three of them which give an approximately
holomorphic generic covering $f_k:X\to \mathbb C\mathbb P^2$ of
degree $N_k=k^2 \omega^2$ branched over a cuspidal Hurwitz curve
$\bar H_k$ (the notion of Hurwitz curves is a generalization of
the notion of plane algebraic curves; see the definition in
section \ref{Hurw}). In algebraic case, if $\deg f\geq 12$, then
$f$ is determined uniquely by $\bar H$ (\cite{Ku}, \cite{Nem}).
Any generic covering $f:X\to \mathbb C\mathbb P^2$ of degree $N$
branched over a cuspidal Hurwitz curve $\bar H$ determines (and is
determined by) its monodromy $\overline{\mu}$, that is, an
epimorphism $\overline{\mu}:\pi_1(\mathbb C\mathbb P^2\setminus
\bar H)\to \mathfrak S_N$ to the symmetric group $\mathfrak S_N$.
Therefore one of the main problems is to investigate the
properties of the fundamental groups of the complements of Hurwitz
curves in order to have a possibility to construct interesting
examples of symplectic 4-manifolds.

In \cite{Ku1}, a class $\mathcal C$ of groups, called $C$-groups,
was defined. This class coincides (see Proposition \ref{cgroup})
with a class consisting of the groups which are
%can be
given by finite presentations of the following form: for some
integer $m$, a function $h:\{1,\dots , m\} \to \mathbb Z$, and a
subset $W=\{ w_{i,j,k}\in \mathbb F_m\, \mid \, 1\leq i,j\leq m,\,
\, 1\leq k\leq h(i,j)\} \subset \mathbb F_m$, where $\mathbb F_m$
is the free group generated by an alphabet $\{x_1,\dots, x_m\}$, a
group $G\in \mathcal C$ possesses the presentation
\begin{equation} \label{zero}
G_W=<x_1,\dots ,x_m \, \mid \, x_i= w_{i,j,k}^{-1}
%(\bar x)
x_jw_{i,j,k}
%(\bar x)
, \, \, w_{i,j,k}\in W\, >.
\end{equation}
Let $\varphi_W: \mathbb F_m\to G_W$ be the canonical epimorphism.
The elements $\varphi_{W}(x_i)\in G$, $1\leq g\leq m$, and the
elements conjugated to them are called the {\it $C$-generators} of
the $C$-group $G$. Let $f:G_1\to G_2$ be a homomorphism of
$C$-groups. It is called a {\it $C$-homomorphism} if the images of
the $C$-generators of $G_1$ under $f$ are  $C$-generators of the
$C$-group $G_2$. We will distinguish the groups $G\in \mathcal C$
up to $C$-isomorphisms. Note that the class $\mathcal C$ contains
the subclasses $\mathcal K$ and $\mathcal L$, respectively, of the
knot and link groups given by Wirtinger presentation. In
\cite{Ku1}, it was proved that the class $\mathcal C$ coincides
with the class of the fundamental groups of the complements of
orientable closed surfaces in the 4-dimensional sphere $S^4$ (with
generalized Wirtinger presentation).

We prove (Theorem \ref{fun}) that the class $\mathcal H=\{ \,
\pi_1(\mathbb C^2\setminus H)\, \}$ of the fundamental groups of
the complements of affine Hurwitz curves $H$ (that is, $H=\bar
H\cap (\mathbb C\mathbb P^2\setminus L_{\infty})$, where $\bar H$
is a Hurwitz curve in $\mathbb C\mathbb P^2$ and $L_{\infty}$ is a
projective line in $\mathbb C\mathbb P^2$ in general position with
respect to $\bar H$) coincides with a subclass of $\mathcal C$
consisting of the $C$-groups $G_W$ with presentation (\ref{zero})
such that the set of the set of words $\{ [x_i,x_1\dots x_m] \}
_{i=1,\dots,m}$ is a subset of $W$. We also describe (Theorem
\ref{tor}) the intersection $\mathcal H\cap \mathcal K$ of the
classes $\mathcal H$ and $\mathcal K$.
%the class of the knot groups.
Besides, we prove (Corollary \ref{cor}) that the class
$\mathcal C$ coincides with a class $\mathcal C_{Br}$ of groups
which possess presentations of the form
$$ G_B=< x_1,\dots, x_m \, \mid  \, x_i^{-1}b(x_i)=\mathbf 1,\, \,
i=1,\dots, m,\, \, b\in B\, > $$ for some $m\in \mathbb N$ and a
finitely generated subgroup $B$ of the braid group $Br_m$, where
$b(x_i)$ is the image of $x_i$ under the standard action of $b\in
Br_m$ on $\mathbb F_m$. The proof of these results is based on a
full twist factorizations formula in the braid group $Br_{2m}$
(see section \ref{semi}). This formula also is applied to prove
the existence of a symplectic 4-manifold $X$ and two isotopic
generic coverings $f_i :X\to \mathbb C\mathbb P^2$, $i=1,2$,
branched, respectively, over cuspidal Hurwitz curves (without
negative nodes) $\bar H_i\subset \mathbb C\mathbb P^2$ having
different braid monodromy factorization types (see Theorem
\ref{last}).

\section{$C$-groups}

Let $\mathbb F$ be the free group generated by an alphabet $\{ x_j
\,\, \mid \, \, j\in \mathbb N\, \, \}$. Below, denote by $\mathbb
F_m$ the subgroup of $\mathbb F$ generated by $x_1,\dots,x_m$, and
$\mathbb F_{m,k}$ the subgroup of $\mathbb F$ generated by
$x_{k+1},\dots,x_{k+m}$. Notation $w(x_1,\dots,x_m)$ means a word
in letters $x_1,\dots,x_m$ and their inverses considered as an
element in the group $\mathbb F_m$.

Let $W$ be a subset of  $\mathbb F_m$ and $N_m(W)$ the normal
closure of $W$ in $\mathbb F_m$. Denote by $\varphi_W: \mathbb
F_m\to \mathbb F_m/N_m(W)=G(W)$ the canonical epimorphism. The set
$$R_{W}=\{ w_{j}(x_1,\dots, x_m)=\mathbf 1 \, \mid \, \, w_{j}\in W\, \}$$
is called a {\it set of relations defining} $G(W)$. We say that
the relations $R_W$ imply a relation $\{ w(x_1,\dots,x_m)=\mathbf
1 \}$ if $w\in N_m(W)$. Sometimes the set $R_W$ will be written in
the form
$$R_W=\{ u_{j}(x_1,\dots,x_m)=v_{j}(x_1,\dots,x_m)\, \mid \,
\, w_{j}(\bar x)=u_{j}(\bar x)v_{j}(\bar x)^{-1}\in W\, \}
$$
for some presentation of the words $w_j\in W$ in the form
$w_{j}=u_{j}v_{j}^{-1}$.

Let $W_1\subset \mathbb F_{m_1}$, $W_2\subset \mathbb F_{m_2}$,
and $f:G(W_1)\to G(W_2)$ be a homomorphism. We say that $f$ is the
{\it canonical homomorphism} if
$$f(\varphi_{W_1}(x_i))=
\varphi_{W_2}(x_i)$$ for all $i\leq \min(m_1,m_2)$. In particular,
if $W_1\subset W_2$ are two subsets in $\mathbb F_m$, then there
is the canonical epimorphism
\begin{equation} \label{epi}
\psi_{W_1,W_2}:G(W_1)\to G(W_2)=
G(W_1)/N_{G(W_1)}(\varphi_{W_1}(W_2\setminus W_1)),
\end{equation}
 where $N_{G(W_1)}(\varphi_{W_1}(W_2\setminus
W_1))$ is the normal closure in $G(W_1)$ of the set
$\varphi_{W_1}(W_2\setminus W_1)$.

We say that two sets of relations
$$R_{W_i}=\{ w_{i,j}=\mathbf 1 \, \mid \, \, w_{i,j}\in W_i\, \} ,\qquad i=1,2,$$
are {\it equivalent} if $N_m(W_1)=N_m(W_2)$.  More generally, for
$k\leq m$, let $W_1\subset \mathbb F_k$ and $W_2\subset \mathbb
F_m$ be two subsets of elements in $\mathbb F_k$ and $ \mathbb
F_m$, respectively. The sets of relations
$$R_{W_1}=\{ w_{1,j}=\mathbf 1 \, \mid \, \, w_{1,j}\in W_1\, \} $$
and
$$R_{W_2}=\{ w_{2,j}=\mathbf 1 \, \mid \, \, w_{2,j}\in W_2\, \} $$
are called {\it equivalent} if  there is the canonical isomorphism
$h: G(W_1)\to G(W_2)$.
%Extend this equivalence to an equivalence relation.

\begin{claim} \label{cl0} For $W_1,W_2\subset \mathbb F_k$ and
$\overline W_1\subset \mathbb F_m$, $m\geq k$, if the sets of
relations $R_{W_1}$ and $R_{\overline W_1}$ are equivalent, then
the sets of relations $R_{W_1\cup W_2}$ and $R_{\overline W_1\cup
W_2}$ are equivalent.
\end{claim}
\proof Straightforward.\qed

\begin{claim} \label{cgr}
Let $w(x_1,\dots,x_m)=x_{j_{1}}^{\varepsilon_{1}}\dots
x_{j_{n}}^{\varepsilon_{n}}$, $\varepsilon_{k}=\pm 1$ for
$k=1,\dots, n$, be a word in the letters $x_1,\dots,x_m$ and their
inverses of the letter length $n\geq 2$. Then for any pair $i_1,i_2\in
\{ 1,\dots,m\}$, the relation $\{ x_{i_1}= w^{-1}x_{i_2}w\}$ in
$\mathbb F_m$ is equivalent to the union of the sets of relations
$$\begin{array}{rl}
R_w= & \{
x_{m+k}=x_{j_k}^{-\varepsilon_k}x_{m+k-1}x_{j_k}^{\varepsilon_k}\,
\mid \, k=2,\dots, n-1\, \} \cup \\
&  \{
x_{m+1}=x_{j_1}^{-\varepsilon_1}x_{i_2}x_{j_1}^{\varepsilon_1},\,
\,
x_{i_1}=x_{j_n}^{-\varepsilon_n}x_{m+n-1}x_{j_n}^{\varepsilon_n}
\}
\end{array}
$$
in $\mathbb F_{m+n-1}$.
\end{claim}
\proof Straightforward.\qed

Let $I$ be any subset in $\{ 1,\dots,N\}^3$. Since for
$(i,j,k)\in \{ 1,\dots,N\}^3$ the relations $x_i=x_j^{-1}x_kx_j$
and $x_k=x_jx_ix_j^{-1}$ are equivalent, Claims \ref{cl0} and
\ref{cgr} imply the following proposition.
\begin{prop} \label{cgroup}
Any $C$-group $G$ is canonically $C$-isomorphic to a $C$-group
$\overline G$ with $C$-presentation of the form
$$\overline G\simeq < \, x_1,\dots,x_N\, \, \mid \, \,
x_i=x_j^{-1}x_kx_j\, \, \mbox{for}\, \, (i,j,k)\in I\, >,
$$
i.e., the definition of $C$-groups given in Introduction coincides
with the definition of $C$-groups given in \cite{Ku1}.
\end{prop}

Note that the relation $x_i=x_j$ is a $C$-relation, since it can
be written in the form $x_i=x_j^{-1}x_jx_j$.
\vspace{0.1cm}
\newline
{\bf Example.} Let $Br_m$ be the braid group on $m$ strings. It is
generated by generators $a_1,\dots ,a_{m-1} $  being subject
to the relations
\begin{equation} \label{cbr}
\begin {array}{ll}
a_ia_{i+1}a_i & =a_{i+1}a_i a_{i+1} \qquad \qquad 1\leq i\leq m-1 ,  \\
a_ia_{k} & =a_{k}a_i  \qquad \qquad \qquad \, \, \mid i-k\mid \,
\geq 2.
\end{array}
\end{equation}
The braid group $Br_m$ with presentation (\ref{cbr}) possesses a
natural  structure of $C$-group. Indeed, let
$$
\begin{array}{l}
R_{W(Br_m)}=\{ \, x_{2i+1}=x_{2i-1}^{-1}x_{2i}x_{2i-1} \, \mid \,
i=1,\dots
,m-1\, \}\cup \\
\{ \, x_{2i-1}=x_{2i}^{-1}x_{2i+1}x_{2i} \, \mid \,
i=1,\dots ,m-1\, \}\cup  \\
\{ \, x_{2i-1}=x_{2j-1}^{-1}x_{2i-1}x_{2j-1}\, \mid \, |i-j|\geq
2,\, \,  1\leq i,j\leq m-1\, \} %\subset \mathbb F_{2m-3},
\end{array}
$$
be a set of $C$-relations in $\mathbb F_{2m-3}$. Put
$G_{Br_m}=G(W(Br_m))$. One can check that the homomorphism $f:
Br_m\to G_{Br_m}$ given by $f(a_i)=\varphi_W(x_{2i-1})$,
$i=1,\dots, m-1$, is an isomorphism.\vspace{0.1cm}

For a set $W=W_f=\{ w_{j}(x_1,\dots,x_m)\in \mathbb F_m \, \mid \,
\, j\in J \}$, where $f:J\to \mathbb F_m$ is a map from a set $J$
and $f(j)$ is denoted by $w_j$ for $j\in J$, the set
$$\mbox{sh}\,(W)=\mbox{sh}_m\,(W)=\{
w_{j}(x_{m+1},\dots,x_{2m}) \, \mid \, \, j\in J\, \}
$$
is called the {\it shift} of $W$ by $m$. Denote by
$$ W_{mod\, m}=\{ x_i^{-1}x_{m+i}\, \, \mid \, \, i=1,\dots, m\, \}$$
and call the sets
\begin{equation} \label{doubl}
d(W)=d_m(W)=W\cup \mbox{sh}\,(W)\cup W_{mod\, m}
\end{equation}
and
\begin{equation} \label{rdoubl}
\overline d(W)=\overline d_m(W)=\mbox{sh}\, (W)\cup W_{mod\, m}
\end{equation}
 in $\mathbb F_{2m}$, respectively, the {\it doubling} and the
 {\it reduced doubling of the set} $W$.

The doubling (respectively, the reduced doubling) of a set $W$ can
be iterated if we put
\begin{equation} \label{iter}
d^n(W)=d_{2^{n-1}m}(d^{n-1}(W))\in \mathbb F_{2^nm}.
\end{equation}

\begin{claim} \label{dclaim}
For any $n\in \mathbb N$ the sets of relations $R_{W}$,
$R_{d^n(W)}$, and $R_{\overline d^n(W)}$ are equivalent.
\end{claim}
\proof Straightforward.\qed

Let $B$ be a subgroup of $\mbox{Aut}\, (\mathbb F_m)$. Denote by
$R_{W(B,\mathbb F_m)}$ the set of relations
$$R_{W(B,\mathbb F_m})=\{ wg(w)^{-1}=\mathbf 1\, \, \mid\, \, w\in \mathbb F_m,\,
\, g\in B\, \}.$$

\begin{claim} \label{claim1} Let $y_1,\dots,y_m\in \mathbb F_m$ be elements
generating $\mathbb F_m$. For a subgroup $B\subset \mbox{Aut}\,
(\mathbb F_m)$, the set of relations $R_{W(B,\mathbb F_m)}$ is
equivalent to
$$R_{W(B,\overline y)}=\{ y_ig(y_i)^{-1}=\mathbf 1\, \, \mid\, \, 1\leq i\leq m,\,
\, g\in B\, \}.$$
\end{claim}

\proof Note that to prove Claim \ref{claim1}, it is sufficient to
show that the condition that elements $w_1g(w_1)^{-1}$ and
$w_2g(w_2)^{-1}$ belong to $N_m(W(B,\overline y))$ for some
$w_1,w_2\in \mathbb F_m$ implies
\[
w_1^{-1}g(w_1^{-1})^{-1}  \in N_m(W(B,\overline y))
\]
and
\[
w_1w_2g(w_1w_2)^{-1}  \in N_m(W(B,\overline y)).
\]

If $w_1g(w_1)^{-1}\in N_m(W(B,\overline y))$, then
$$(w_1g(w_1)^{-1})^{-1}=g(w_1)w_1^{-1}\in N_m(W(B,\overline y))$$
and therefore
$$w_1^{-1}g(w_1^{-1})^{-1}=w_1^{-1}(g(w_1)w_1^{-1})w_1\in N_m(W(B,\overline y)).$$
If $w_1g(w_1)^{-1},w_2g(w_2)^{-1}\in N_m(W(B,\overline y))$, then
\[
\begin{array}{l}
w_1w_2g(w_1w_2)^{-1}=w_1w_2g(w_2)^{-1}g(w_1)^{-1}= \\
w_1(w_2g(w_2)^{-1})w_1^{-1}(w_1g(w_1)^{-1}) \in N_m(W(B,\overline
y)).
\end{array}
\]
%Claim \ref{claim1} is proved.
\qed

\begin{claim} \label{claim2} If $B\subset
\mbox{Aut}\, (\mathbb F_m)$ is generated by elements $b_1,\dots,
b_n$ and elements $y_1,\dots,y_m$  generate $\mathbb F_m$, then
$R_{W(B,\mathbb F_m)}$ is equivalent to
$$R_{W(\overline b,\overline y)}=
\{ y_ib_j(y_i)^{-1}=\mathbf 1\, \, \mid\, \, 1\leq i\leq m,\,
\, 1\leq j\leq n\, \}.$$
\end{claim}
\proof It follows from the proof of Claim \ref{claim1} that to
prove Claim \ref{claim2}, it is sufficient to show that the
condition that the elements $wg_1(w)^{-1}$ and $wg_2(w)^{-1}$
belong to $N_m(W(\overline b,\overline y))$ for some $g_1,g_2\in
B$ and any $w\in \mathbb F_m$ implies $w(g^{-1}_1(w))^{-1}\in
N_m(W(\overline b,\overline y))$ and $w(g_1g_2(w))^{-1} \in
N_m(W(\overline b,\overline y))$.

If $wg_1(w)^{-1}\in N_m(W(\overline b,\overline y))$ for all $w\in
\mathbb F_m$, then
\[
\begin{array}{l}
(w(g_1^{-1}(w))^{-1})^{-1}= (g_1^{-1}(w))w^{-1}= \\
(g_1^{-1}(w))(g_1(g_1^{-1}(w)))^{-1}\in N_m(W(\overline
b,\overline y)).
\end{array}
\]

If $w(g_1(w))^{-1},w(g_2(w))^{-1} \in N_m(W(\overline b,\overline
y))$ for all $w\in \mathbb F_m$, then
\[
\begin{array}{l}
wg_1g_2(w)^{-1}=wg_1(g_2(w))^{-1}= \\
(wg_2(w)^{-1})(g_2(w)g_1(g_2(w))^{-1}) \in N_m(W(\overline
b,\overline y)).
\end{array}
\]
%Claim \ref{claim2} is proved.
\qed

The braid group $Br_m$ on $m$ strings acts on $\mathbb F_m$.
Below, we fix a set $\{ a_1,\dots ,a_{m-1} \} $ of so called {\it
standard generators}, i.e., generators of $Br_m$ being subject to
the relations
$$
\begin {array}{ll}
a_ia_{i+1}a_i & =a_{i+1}a_i a_{i+1} \qquad \qquad 1\leq i\leq n-1 ,  \\
a_ia_{k} & =a_{k}a_i  \qquad \qquad \qquad \, \, \mid i-k\mid \,
\geq 2
\end{array}
$$
and acting on $\mathbb F_m$ as follows:
\[
\begin{array}{ll}
a_i(x_j)=x_j\quad &\mbox{if}\, \, i\neq j,j+1; \\
a_i(x_i)=x_ix_{i+1}x_i^{-1}; &  \\
a_i(x_{i+1})=x_{i}. & \end{array} \]
Therefore one can associate a
set of relations $R_{W(B,\mathbb F_m)}$ to any subgroup $B\subset
Br_m$. As in the case of the free groups, for $k\leq m $ the group $Br_k$ is identified
with the subgroup of $Br_m$ generated by the first $k-1$ standard
generators $a_1,\dots,a_{k-1}$, and we denote by $B_{k,i}$,
$k+i\leq m$, the subgroup of $Br_m$ generated by $a_{i+1},\dots,
a_{i+k-1}$. Let $\Delta _{k,i}$  be the Garside element in the
braid group $B_{k,i}$:
$$
\Delta_{k,i}= (a_{i+1}\dots
a_{i+k-1})\dots(a_{i+1}a_{i+2}a_{i+3})(a_{i+1}a_{i+2})a_{i+1}.
$$
We have
\begin{equation}
\label{delta} \Delta_{k,i}(x_{i+j})=(x_{i+1}\dots
x_{i+k-j})x_{i+k-j+1}(x_{i+1}\dots x_{i+k-j})^{-1}
\end{equation}
for $i+1\leq j\leq i+k$ (we put here $x_{i+1}\dots
x_{i+k-j}=\mathbf 1$ if $j=i+k$) and
$$\Delta_k(x_j)=x_j$$
for $j\leq i$ and $j>i+k$. The Garside element in $Br_{k}$ will be
denoted by $\Delta_{k}$.

Let $\mbox{sh}=\mbox{sh}_m: Br_m\to Br_{2m}$ be the homomorphism
given by $\mbox{sh}(a_j)=a_{m+j}$ for $j=1,\dots,m-1$. For any
subgroup $B$ of $Br_m$, its image $\mbox{sh}(B)$ will be called
the {\it shift} of the group $B$ by $m$.

Below, we will identify $\mathbb F_k$ with the fundamental group
$\pi_1(D_k\setminus \{p_1,\dots,p_k \},p_0)$, where $D_k=\{ z\in
\mathbb C\, \mid \, |z|\leq k+1\, \}$, $p_0=-i(k+1)$, $p_j=-j$ for
$j=1,\dots, k$, and for each $j$ the element $x_j$ will be
identified with the loop consisting of the segment $\gamma_j=\{ z\in
\mathbb C\, \mid \, z=-i(k+1)- t(j+i(\varepsilon -k-1)), 0\leq
t\leq 1\, \}$ for some $0<\varepsilon <<1$, the circuit around the
circle $\{ |z+j|=\varepsilon \}$ in the counterclockwise direction and
the return along $l_j$ (see Fig. 1). For so chosen base
$x_1,\dots, x_k$, the standard generators $a_j$ of $Br_k$,
$j=1,\dots,k-1$, are identified with half-twists defined by the
segments $[-j-1,-j]=\{ z=-(j+1)+t\, \mid\, 0\leq t\leq 1\, \}$
(see details, for example, in \cite{Moi-Te2}). The element
$l_k=x_1\dots x_k$ coincides in $\pi_1(D_k\setminus \{p_1,\dots,p_k
\},p_0)$ with circuit around the boundary $\partial D_k$ of $D_k$
in the counterclockwise direction and it is fixed under the action
of $Br_k$.
\begin{picture}(400,310)

\put(180,140){\circle*{3}}
\put(160,140){\circle*{3}}
%\put(140,140){\circle*{3}}
%\put(120,140){\circle*{3}}

\put(100,140){\circle*{3}}

%\put(80,140){\circle*{3}}

\put(60,140){\circle*{3}}

\put(182,15){$p_0$}
%\put(180,280){\circle*{3}}
%\put(182,283){$k+1$}
\put(180,13){\circle*{3}}
%\put(320,140){\circle*{3}}
%\put(310,147){$k+1$}
\put(180,0){\vector(0,1){290}} \put(40,140){\vector(1,0){290}}
\put(160,140){\circle{10}} \put(160,145){\vector(-1,0){1}}
\put(147,147){$-1$} \qbezier(180,13)(170,75)(160,137)
\put(170,95){$x_{1}$} \put(100,140){\circle{10}}
\put(96,147){$-i$} \qbezier(180,13)(140,80)(100,136)
\put(130,95){$x_{i}$} \put(60,140){\circle{10}} \put(54,147){$-k$}
\qbezier(180,13)(120,75)(60,136) \put(103,95){$x_{k}$}
\end{picture}

\begin{center} Fig. 1 \end{center}

We say that an ordered set $\{ y_1,\dots, y_m\}$ consisting of elements of
$\mathbb F_m$ is a {\it good geometric base} of $\mathbb F_m$ if
there is $b\in Br_m$ such that $y_j=b(x_j)$ for $j=1,\dots,m$. An
ordered set $\{ y_1,\dots, y_k\}$, $k\leq m$, of elements of
$\mathbb F_m$ is called {\it good} if the set $\{ y_1,\dots,
y_k\}$ can be extended to a good geometric base $\{ y_1,\dots,
y_m\}$ of $\mathbb F_m$. It is well-known that $\{ y_1,\dots,
y_m\}$ is a good geometric base of $\mathbb F_m$ if and only if
each $y_j$ is conjugated to some $x_{k(j)}$ in $\mathbb F_m$ and
$y_1\dots y_m=l_m$.

Consider an element
$w(x_1,\dots,x_m)=x_{j_{1}}^{\varepsilon_{1}}\dots
x_{j_{n}}^{\varepsilon_{n}}$, $\varepsilon_{j}=\pm 1$, in $\mathbb
F_m$ of the letter length $n$. For $k>n$ and $p\geq 1$,  the
element
$$\mbox{m}_{pmk,k,m}(w)
=x_{(pk-\varepsilon_{1}(k-1))m+j_{1}}^{\varepsilon_{1}}\dots
x_{(pk-\varepsilon_{n}(k-n))m+j_{n}}^{\varepsilon_{n}}\in \mathbb
F_{m(p+1)k}
$$
is called the {\it moving apart} of $w$ by $m$ with center at
$pmk$.

Below, we will use the following claim which is obvious from
geometric point of view. %(see Fig. 2).
\begin{claim}\label{move}
Let $w_k(x_1,\dots,x_m)=x_{j_{k,1}}^{\varepsilon_{k,1}}\dots
x_{j_{k,n_k}}^{\varepsilon_{k,n_k}}$, $\varepsilon_{k,j}=\pm 1$,
$k=1,2$, be two words in the letters $x_1,\dots,x_m$ and their
inverses of the letter length $n_k$ less than $n$. Then for each
pair $i_1,i_2$ such that $1\leq i_1,i_2\leq m$, the pair
$(y_1=\overline w_1x_{nm+i_1}\overline w_1^{-1},y_2=\overline
w_2x_{3nm+i_2}\overline w_2^{-1})$  is good in $\mathbb F_{4mn}$,
where
$$\overline w_1=\mbox{\rm m}_{mn,n,m}(w_1)
=x_{(n-\varepsilon_{1,1}(n-1))m+j_{1,1}}^{\varepsilon_{1,1}}\dots
x_{(n-\varepsilon_{1,n_1}(n-n_1))m+j_{1,n_1}}^{\varepsilon_{1,n_1}}$$
and
$$\overline w_2=\mbox{\rm m}_{3mn,n,m}(w_1)=
x_{(3n-\varepsilon_{2,1}(n-1))m+j_{2,1}}^{\varepsilon_{2,1}}\dots
x_{(3n-\varepsilon_{2,n_2}(n-n_2))m+j_{2,n_2}}^{\varepsilon_{2,n_2}}.$$
\end{claim}

\proof The elements $y_1$ and $y_2$ are represented by loops of the form
drawn in Fig. 2. Obviously, such pair $(y_1,y_2)$ is good. \qed

\begin{picture}(500,170)

\put(262,140){\circle*{3}} \put(233,140){$\dots$}
\put(220,140){\circle*{3}} \put(220,140){\circle{10}}
\put(195,140){$\dots$} \put(183,140){\circle*{3}}
%\put(160,140){\circle*{3}}
\put(120,140){\circle*{3}} \put(99,140){$\dots$}
\put(87,140){\circle*{3}} \put(63,140){$\dots$}
\put(54,140){\circle*{3}} \put(54,140){\circle{10}}
\put(32,140){$\dots$} \put(18,140){\circle*{3}}
%\put(20,140){\circle{10}}
\qbezier(12,140)(12,40)(350,10) \qbezier(12,140)(12,160)(25,140)
\qbezier(54,135)(69,100)(80,140) \qbezier(80,140)(87,160)(94,140)
\qbezier(94,140)(103,120)(114,140)
\qbezier(114,140)(125,160)(126,140)
\qbezier(126,140)(110,40)(25,140) \put(130,140){$\dots$}
\put(149,140){\circle*{3}} \put(140,133){$\mbox{}_{-2nm}$}
%\qbezier(270,10)(30,60)(30,140)
\put(156,140){$\dots$} \qbezier(175,140)(175,160)(190,140)
\qbezier(220,135)(240,120)(252,140)

\put(350,10){\circle*{3}} \put(354,10){$p_0$}
%\put(330,10){\line(0,1){130}}
\qbezier(350,10)(330,40)(310,140)
\qbezier(310,140)(310,160)(292,140)
\qbezier(292,140)(175,-60)(175,140)
\qbezier(190,140)(270,70)(270,140)
\qbezier(270,140)(270,155)(252,140) \put(276,140){$\dots$}
\put(302,140){\circle*{3}} \put(317,140){$\dots$}
\put(337,140){\circle*{3}} \put(331,133){$\mbox{}_{-1}$}
\put(-5,140){$\dots$} \put(-12,140){\circle*{3}}
\put(-20,133){$\mbox{}_{-4nm}$} \put(310,80){$y_{1}$}
\put(60,80){$y_2$}

\end{picture}

\begin{center} Fig. 2 \end{center}

Let $h:\{ 1,\dots, m\} ^2\to \mathbb Z$ be a non-negative
integer-valued function. Put
$$I_h=\{ (i,j,k)\in \mathbb Z^3 \,
\mid \, 1\leq i,j\leq m,\, \,  1\leq k\leq h(i,j)\, \}$$
and let $f:I_h \to \mathbb F_m^2$ be a map to
$\mathbb F_m^2=\mathbb F_m\times \mathbb F_m$.

\begin{claim}\label{main1} For any map $f:I_h \to \mathbb F_m^2$,
there are $M\in
\mathbb N$, a finite set $I_{h,f}\subset \{1,\dots,M\}^2$ and a
map $F:I_{h,f}\to \mathbb F_M^2$ such that
\begin{itemize}
\item[(i)] the pair $F(i,j)=(y_i,y_j)$ is a good pair in $\mathbb F_M$
for each $(i,j)\in I_{h,f}$;
\item[(ii)] the sets of relations
$$R_{W_{f}}=
\{\, x_i^{-1}w_{i,j,k}x_jw_{i,j,k}^{-1}=\mathbf 1 \, \mid \,
(i,j,k)\in I_h, \,\, w_{i,j,k}=f(i,j,k)\, \}$$ and
$$R_{W_F}=\{ \, y_i^{-1}y_j=\mathbf 1\, \mid \, (y_i,y_j)\in
\mbox{\rm Im}\, F\, \}$$ are equivalent.
\end{itemize}
\end{claim}
\proof By Claim \ref{move}, there is an integer $N$ such that the
pairs $$(x_i,
\mbox{m}_{Nm,N,m}(w_{i,j,k})x_{Nm+j}\mbox{m}_{Nm,N,m}(w_{i,j,k})^{-1})$$
are good for all $(i,j,k)\in I_h$.  Obviously, the set of
relations  $R_{W_{f}}$ in $\mathbb F_m$ is equivalent to the union
of two sets of relations $R_{\mbox{\rm m}(W_f)}$ and $R_{mod\, m}$
in $\mathbb F_{2Nm}$, where
$$R_{\mbox{\rm m}(W_f)}=\{x_i= \mbox{m}_{Nm,N,m}
(w_{i,j,k})x_{Nm+j}\mbox{m}_{Nm,N,m}(w_{i,j,k})^{-1}\,
\mid \, \, (i,j,k)\in I_h\, \}$$ and
$$ R_{mod\, m}=\{ x_i=x_j \, \mid \,
i\equiv j \mod\, m)\} .$$ For each $i\leq m$ and each $k> 0$ the
pair $(x_i,x_{i+mk})$ is good which completes the proof of Claim
\ref{main1} if we put $M=2Nm$ and $R_{W_F}=R_{\mbox{\rm
m}(W_f)}\cup R_{mod\, m}$.\qed

By Claim \ref{claim2}, for a subgroup $B$ of $Br_m$ generated by
$b_1,\dots, b_n$, the group $G(W(B,\mathbb F_m))$ is canonically
isomorphic to
\begin{equation} \label{conj2}
\begin{array}{l}
G(W(\overline b,\overline x))= \\
<x_1, \dots , x_m \, \mid \,
x_i^{-1}b_j(x_i)=\mathbf 1, \, \, i=1,\dots, m,\, \, \, j=1,\dots,
n
>.
\end{array}
\end{equation}

\begin{claim} \label{newclaim}
Presentation (\ref{conj2}) defines on $G(W(\overline b,\overline
x))$ a structure of a $C$-group.
\end{claim}
\proof We have
\begin{equation} \label{action} b(x_i)=w_{b,i}(\bar
x)^{-1}x_jw_{b,i}(\bar x)
\end{equation}
for some $j=\sigma_b(i)$ and $w_{b,i}(\bar x)\in \mathbb F_m$,
where $\sigma_b\in \mathfrak S_m$ is the image of $b\in Br_m$
under the canonical epimorphism $\sigma :Br_m\to \mathfrak S_m$
with $\ker \sigma =P_m$ (we denote by $\mathfrak S_m$ the
symmetric group acting on the set $\{ 1,\dots, m\}$ and by $P_m$
the subgroup of pure braids). Therefore, the relations
$x_i^{-1}b_j(x_i)=\mathbf 1$ are equivalent to the $C$-relations
$x_i =w_{b_j,i}(\bar x)^{-1}x_jw_{b_j,i}(\bar x)$. \qed

The following lemma is well known.
\begin{lem}\label{Ak} Let
$b_k=ga_j^{k+1}g^{-1}\in Br_m$. Then the set of relations
$$\{ x_i^{-1}b_k(x_i)=\mathbf 1\, \mid  \, \, i=1,\dots, m\}$$
in $\mathbb F_m$ is equivalent to one single relation:
\begin{itemize}
\item[(0)]
$g(x_j)=g(x_{j+1})$ if $k=0$;
\item[(1)] the commutant
$[g(x_j),g(x_{j+1})]=\mathbf 1$ if $k=1$ or $-3$;
\item[(2)] $g(x_j)g(x_{j+1})g(x_{j})=g(x_{j+1})g(x_j)g(x_{j+1})$ if $k=2$.
\end{itemize}
\end{lem}

\proof Consider a free base $\{ y_1,\dots y_{m}\}$ of the group
$\mathbb F_{m}$, where $y_i=g(x_i)$ for $i=1,\dots,m$. By Claim
\ref{claim2}, the sets of relations $\{ x_i^{-1}b_k(x_i)=\mathbf
1\, \mid  \, \, i=1,\dots, m\}$ and $\{ y_i^{-1}b_k(y_i)=\mathbf
1\, \mid  \, \, i=1,\dots, m\}$ are equivalent for each $k$.

We have $$b_k(y_i)=g(x_i)=ga_j^{k+1}g^{-1}(g(x_i))=g(x_i)=y_i$$ if
$i\neq j,j+1$,
\[
\begin{array}{ll}
b_k(y_j) & =ga_j^{k+1}(x_j)= \\
 & \left\{
\begin{array}{ll}
g(x_jx_{j+1}x_j^{-1})\, &
\mbox{if}\, \, k=0; \\
g((x_jx_{j+1})x_j(x_{j}x_{j+1})^{-1})\, &
\mbox{if}\, \, k=1; \\
(g((x_jx_{j+1}x_j)x_{j+1}(x_jx_{j+1}x_j)^{-1})\, & \mbox{if}\, \,
k=2,
\end{array} \right.
\end{array}
\]
and
\[
\begin{array}{ll}
b_k(y_{j+1}) & =ga_j^{k+1}(x_{j+1})= \\
 & \left\{
\begin{array}{ll}
g(x_j)\, &
\mbox{if}\, \, k=0; \\
g(x_jx_{j+1}x_j^{-1})\, &
\mbox{if}\, \, k=1; \\
g((x_jx_{j+1})x_j(x_{j+1}x_j)^{-1})\, & \mbox{if}\, \, k=2.
\end{array} \right.
\end{array}
\]
Therefore
\[
b_k(y_j) =
  \left\{
\begin{array}{ll}
y_jy_{j+1}y_j^{-1}\, &
\mbox{if}\, \, k=0; \\
(y_jy_{j+1})y_j(y_{j}y_{j+1})^{-1}\, &
\mbox{if}\, \, k=1; \\
(y_jy_{j+1}y_j)y_{j+1}(y_jy_{j+1}y_j)^{-1}\, & \mbox{if}\, \, k=2
\end{array} \right.
\]
and
\[
b_k(y_{j+1})=  \left\{
\begin{array}{ll}
y_j\, &
\mbox{if}\, \, k=0; \\
y_jy_{j+1}y_j^{-1}\, &
\mbox{if}\, \, k=1; \\
(y_jy_{j+1})y_j(y_{j+1}y_j)^{-1}\, & \mbox{if}\, \, k=2.
\end{array} \right.
\]
Thus, for each $k$ the set of relations $\{
y_i^{-1}b_k(y_i)=\mathbf 1\, \mid \, \, i=1,\dots, m\}$ is reduced
to one single relation
\[
\begin{array}{rl}
y_jy_{j+1}^{-1}=\mathbf 1 \, \, &
\mbox{if}\, \, k=0; \\
y_jy_{j+1}y_j^{-1}y_{j+1}^{-1}=\mathbf 1 \, \, &
\mbox{if}\, \, k=1; \\
(y_jy_{j+1}y_j)(y_{j+1}y_jy_{j+1})^{-1}=\mathbf 1 \, \, & \mbox{if}\, \, k=2.
\end{array}
\]
% Lemma \ref{Ak} is proved.
\qed

\begin{claim} \label{pair}
Let $y_1,y_2$ be a good pair in $\mathbb F_m$. Then there is an
element $b=ga_1g^{-1}\in Br_m$ such that the set of relations
$R_{W(<b>,\mathbb F_m)}$ is equivalent to the single relation
$$R_b=\{ y_1^{-1}y_2=\mathbf 1 \}.$$
\end{claim}
\proof  Since $y_1,y_2$ is a good pair, it can be included in a
good geometric base $\{ y_1,\dots,y_m \}$ of $\mathbb F_m$.
Therefore there is an element $g\in Br_m$ such that $g(x_i)=y_i$
for $i=1,\dots,m$.  By Lemma \ref{Ak}, $R_{W(<b>,\mathbb F_m)}$ is
equivalent to $R_b$. %Claim \ref{pair} is proved.
\qed

Consider the braid group $Br_{2m}$. Put $a_{m,i}=a_ia_{2m-i}$,
\begin{equation} \label{barc} \bar c_{m,i}=(a_{m,i}a_{m,i+1}\dots a_{m,m-1})
a_m(a_{m,i}a_{m,i+1}\dots a_{m,m-1})^{-1}
\end{equation}
for $1\leq i\leq m-1$, $\bar c_{m,m}=a_m$, and
\begin{equation}
\label{ci} c_{m,i}=\Delta_{m,m}^{-1}\bar c_{m,i}\Delta_{m,m}.
\end{equation}
\begin{claim} \label{shift}
The set of relations $R_{W(<c_{m,i}>,\mathbb F_{2m})}$ is
equivalent to the single relation
$$R_{c_{m,i}}=\{ \, x_{m+i}^{-1}x_{i}=\mathbf 1\, \} .$$
\end{claim}
\proof Using (\ref{delta}), one can check that
$$\Delta_{m,m}^{-1}(a_{m,i}a_{m,i+1}\dots a_{m,m-1})(x_{m})=x_i$$ and
$$\Delta_{m,m}^{-1}(a_{m,i}a_{m,i+1}\dots a_{m,m-1})(x_{m+1})=x_{m+i}.$$
Therefore Claim \ref{shift} follows from (\ref{barc}), (\ref{ci}),
and Claim \ref{Ak}.\qed

For any subgroup $B$ of $Br_m$, denote by $d(B)=d_m(B)$
(respectively, $\overline d(B)=\overline d_m(B)$) the subgroup of
$Br_{2m}$ generated by the elements of $B$, $\mbox{sh}_m(B)$, and
the elements $c_{m,1},\dots,c_{m,m}$ (respectively, generated by
the elements of $B$ and the elements $c_{m,1},\dots,c_{m,m}$), and
call it the {\it doubling} (respectively, the {\it reduced
doubling}) of the group $B$.

The doubling (respectively, the reduced doubling) of a group $B$
can be iterated if we put
\begin{equation} \label{iter2}
d^n(B)=d_{2^{n-1}m}(d^{n-1}(B))\subset Br_{2^nm}.
\end{equation}

\begin{claim} \label{dclaim2}
For any $n\in \mathbb N$ the following sets of relations
\[R_{W(B,\mathbb F_m)}, \, R_{d^n(W(B,\mathbb F_m))}, \, R_{\overline
d^n(W(B,\mathbb F_m))}, \, R_{W(d^n(B),\mathbb F_{2^nm})},  \,
\mbox{and} \, R_{W(\overline d^n(B),\mathbb F_{2^nm})}\] are
equivalent.
\end{claim}
\proof It follows from Claims \ref{dclaim} and \ref{shift}.\qed

We say that a $C$-group $G$ belongs to a subclass $\mathcal
C_{Br}$ if it is $C$-isomorphic to a group given by the
presentation
$$G=< x_1,\dots,x_m\, \mid \, R_{W(B,\mathbb F_m)}\, >$$
for some $m$ and some finitely generated subgroup $B$ of $Br_m$.

\begin{thm} For any $C$-group $G$ there are $M_G\in \mathbb Z$ and a subgroup $B_G$ of
$Br_M$ generated by a finite set $\{ b_1,\dots, b_{n_G}\}$ of
elements conjugated in $Br_M$ to the standard generator $a_1$ such
that the group $G$ is $C$-isomorphic to $G(W(B,\mathbb F_M))$.
\end{thm}
\proof It follows from Claims \ref{dclaim}, \ref{move}
\ref{main1}, \ref{pair}, and \ref{shift}.\qed

\begin{cor} \label{cor}
$\mathcal C_{Br}=\mathcal C$.
\end{cor}

Recall also the following theorem.
\begin{thm} {\rm (\cite{Ar})}
\begin{itemize}
\item[(i)]
A $C$-group $G\in \mathcal L$ if and only if there are an integer
$m$ and an element $b\in Br_m$ such that $G$ is $C$-isomorphic to
\begin{equation} \label{link}
G(W(<b>,\bar x))=<\, x_1,\dots,x_m\, \mid \, x_i=b(x_i)\, \,
\mbox{for}\, \, i=1,\dots,m\, >.
\end{equation}
\item[(ii)]
A $C$-group $G$ with presentation {\rm (\ref{link})} belongs to
$\mathcal K$ if and only if the permutation $\sigma_b\in \mathcal
S_m$ consists of one cycle of the length $m$, i.e., $\sigma_b$
acts transitively on the set $\{1,\dots,m\}$.
\end{itemize}
\end{thm}

\section{Full twist factorizations}
It is well-known that the element (called the {\it full twist})
$$ \Delta_m^2=(a_1\dots a_{m-1})^{m}$$
is the generator of the center of $Br_m$. We have
\begin{equation}
\label{delta2} \Delta_m^2(x_i)=(x_1\dots x_m)x_{i}(x_1\dots
x_m)^{-1}.
\end{equation}
In particular, the full twist $\Delta_m^2$ leaves fixed the
element $l_m=x_1\dots x_m$.

Consider the braid group $Br_{2m}$. Put
$$\bar r_{m}=\bar c_{m,m}\bar c_{m,m-1}\dots \bar c_{m,1}$$
and
$$r_{m}=c_{m,m}c_{m,m-1}\dots c_{m,1},$$
where $\bar c_{m,i}$ and $c_{m,i}$ were defined by means of
(\ref{barc}) and (\ref{ci}).

\begin{lem} \label{Delta} In the braid group $Br_{2m}$, we have
\begin{equation}
\Delta_{2m}=\Delta^2_{m,0}\Delta^2_{m,m}\bar r_{m}.
\end{equation}
\end{lem}

\proof  Identify $Br_{2(m-1)}$ with subgroup $B_{2(m-1),1}\subset
Br_{2m}$. By induction on $m$, we have
$$\bar r_{m}=\bar c_{m,m}\bar c_{m,m-1}\dots \bar c_{m,2}\bar c_{m,1}=
\Delta_{(m-1),m}^{-2}\Delta_{m-1,1}^{-2}\Delta_{2(m-1),1}\bar
c_{m,1}.$$ Therefore to prove Lemma \ref{Delta}, it is sufficient
to show that the actions on
$\mathbb F_{2m}$ of the elements $\Delta_{2m}$ and
$$b=\Delta^2_{m,0}\Delta^2_{m,m}\Delta_{m-1,m}^{-2}
\Delta_{m-1,1}^{-2}\Delta_{2(m-1),1}\bar c_{m,1}$$
coincide.

We fix the base $\{ x_1,\dots,x_{2m}\}$ of $\mathbb F_{2m}$ which is
drawn in Fig. 1. Then the action of $\bar c_{m,1}$ on $\mathbb F_{2m}$
is induced by the
half-twist defined by the path $\gamma_{m,1}$ which is drawn in
Fig. 3.

\begin{picture}(500,200)

\put(280,140){\circle*{3}} \put(274,147){$\mbox{}_{-1}$}
\put(260,140){\circle*{3}} \put(240,140){\circle*{3}}
\put(220,140){\circle*{3}} \put(200,140){\circle*{3}}
\put(180,140){\circle*{3}} \put(173,147){$\mbox{}_{-m}$}
\put(140,140){\circle*{3}} \put(123,131){$\mbox{}_{-m-1}$}
\put(120,140){\circle*{3}} \put(100,140){\circle*{3}}
\put(80,140){\circle*{3}} \put(60,140){\circle*{3}}
\put(40,140){\circle*{3}}\put(25,131){$\mbox{}_{-2m}$}

\qbezier(40,140)(100,240)(160,140)
\qbezier(160,140)(220,40)(280,140) \put(220,80){$\gamma_{m,1}$}
%\put(130,100){$\mbox{Fig. 1}$}
\end{picture}

\begin{center} Fig. 3\end{center}
Therefore the element $y_1=\bar c_{m,1}(x_1)$ is represented by
the loop drawn in Fig. 4.

\begin{picture}(500,200)

\put(260,140){\circle*{3}} \put(254,131){$\mbox{}_{-1}$}
\put(234,147){$\mbox{}_{-2}$} \put(240,140){\circle*{3}}
%\put(220,140){\circle*{3}} \put(200,140){\circle*{3}}
\put(195,140){$\dots$} \put(160,140){\circle*{3}}
\put(153,147){$\mbox{}_{-m}$}
%\put(140,140){\circle*{3}}
\put(120,140){\circle*{3}} \put(103,131){$\mbox{}_{-m-1}$}
%\put(100,140){\circle*{3}} \put(80,140){\circle*{3}}
%\put(60,140){\circle*{3}}
\put(65,140){$\dots$}
\put(20,140){\circle*{3}}\put(5,131){$\mbox{}_{-2m}$}
\put(20,140){\circle{10}} \qbezier(250,140)(270,180)(270,140)
\put(270,10){\circle*{3}} \put(274,10){$p_0$}
\put(270,10){\line(0,1){130}} \qbezier(20,145)(80,240)(140,140)
\qbezier(140,140)(200,40)(250,140) \put(200,80){$y_{1}$}
\end{picture}
\begin{center}
Fig. 4
\end{center}
The element $u_1=\Delta_{2(m-1),1}(y_1)$ is represented by the
loop drawn in Fig. 5.

\begin{picture}(500,220)
\put(260,140){\circle*{3}} \put(254,131){$\mbox{}_{-1}$}
%\put(234,147){$\mbox{}_{-2}$}
%\put(240,140){\circle*{3}} \put(220,140){\circle*{3}}
%\put(200,140){\circle*{3}}
\put(205,140){$\dots$}
\put(160,140){\circle*{3}}
\put(153,131){$\mbox{}_{-m}$}
%\put(140,140){\circle*{3}}
\put(120,140){\circle*{3}} \put(103,149){$\mbox{}_{-m-1}$}
%\put(100,140){\circle*{3}} \put(80,140){\circle*{3}}
%\put(60,140){\circle*{3}}
\put(65,140){$\dots$}
\put(20,140){\circle*{3}}\put(7,149){$\mbox{}_{-2m}$}
\put(20,140){\circle{10}}
%\qbezier(250,140)(260,160)(270,140)
\put(270,10){\circle*{3}} \put(274,10){$p_0$}
\put(270,10){\line(0,1){130}} \qbezier(20,135)(80,40)(140,140)
\qbezier(140,140)(250,280)(270,140) \put(257,80){$u_{1}$}
\end{picture}
\begin{center}
Fig. 5
\end{center}
It is easy to see that $\Delta_{m-1,m}^{-2}
\Delta_{m-1,1}^{-2}(u_1)=u_1$ and the element
$v_1=\Delta^2_{m,0}\Delta^2_{m,m}(u_1)=b(x_1)$ is represented by
the loop drawn in Fig. 6.
\begin{picture}(500,220)
\put(153,131){$\mbox{}_{-m}$}
 \put(103,149){$\mbox{}_{-m-1}$}
\put(260,140){\circle*{3}} \put(254,131){$\mbox{}_{-1}$}
%\put(234,147){$\mbox{}_{-2}$}
%\put(240,140){\circle*{3}} \put(220,140){\circle*{3}}
%\put(200,140){\circle*{3}}
\put(205,140){$\dots$}
\put(160,140){\circle*{3}} %\put(153,131){$\mbox{}_{-m}$}
%\put(140,140){\circle*{3}}
\put(120,140){\circle*{3}} %\put(103,149){$\mbox{}_{-m-1}$}
%\put(100,140){\circle*{3}} \put(80,140){\circle*{3}}
%\put(60,140){\circle*{3}}
\put(65,140){$\dots$}
\put(20,140){\circle*{3}}\put(7,131){$\mbox{}_{-2m}$}
\put(20,140){\circle{10}}
%\qbezier(250,140)(260,160)(270,140)
\put(270,10){\circle*{3}} \put(274,10){$p_0$}
\put(270,10){\line(0,1){130}} \qbezier(20,145)(240,270)(270,140)
%\qbezier(140,140)(200,240)(270,140)
\put(257,80){$v_{1}$}
\end{picture}
\begin{center}
Fig. 6
\end{center}
Therefore, by (\ref{delta}),
\begin{equation} \label{first}
b(x_1)=\Delta_{2m}(x_1).
\end{equation}

The element $y_i=\bar c_{m,1}(x_i)$, $i=2,\dots, m$, is
represented by the loop drawn in Fig. 7.

\begin{picture}(500,220)
\put(260,140){\circle*{3}} \put(254,131){$\mbox{}_{-1}$}
\put(234,147){$\mbox{}_{-2}$} \put(240,140){\circle*{3}}
\put(215,140){$\dots$} \put(200,140){\circle*{3}}
\put(160,140){\circle*{3}} \put(175,140){$\dots$}
\put(153,147){$\mbox{}_{-m}$}
%\put(140,140){\circle*{3}}
\put(120,140){\circle*{3}} \put(103,131){$\mbox{}_{-m-1}$}
%\put(100,140){\circle*{3}} \put(80,140){\circle*{3}}
\put(75,140){$\dots$} \put(55,140){\circle*{3}}
\put(20,140){\circle*{3}}\put(44,132){$\mbox{}_{1-2m}$}
 \qbezier(250,140)(270,170)(270,140)
\put(270,10){\circle*{3}} \put(274,10){$p_0$}
\put(270,10){\line(0,1){130}} \qbezier(10,140)(20,260)(150,140)
\qbezier(150,140)(200,70)(200,135) \put(200,140){\circle{10}}
\put(200,131){$\mbox{}_{-i}$} \qbezier(28,140)(80,180)(130,140)
\qbezier(10,140)(10,120)(28,140)
\qbezier(130,140)(200,40)(250,140) \put(273,80){$y_{i}$}
\end{picture}

\begin{center}
Fig. 7
\end{center}
The element $u_i=\Delta_{2(m-1),1}(y_i)$, $i=2,\dots, m$, is
represented by the loop drawn in Fig. 8.

\begin{picture}(500,220)
\put(260,140){\circle*{3}} \put(254,131){$\mbox{}_{-1}$}
%\put(234,147){$\mbox{}_{-2}$}
%\put(240,140){\circle*{3}}\put(220,140){\circle*{3}}
\put(205,140){$\dots$}
%\put(180,140){\circle*{3}}
\put(160,140){\circle*{3}}
\put(153,131){$\mbox{}_{-m}$}
%\put(140,140){\circle*{3}}
\put(120,140){\circle*{3}} \put(103,147){$\mbox{}_{-m-1}$}
\put(95,140){$\dots$} \put(80,140){\circle*{3}}
\put(55,140){$\dots$} \put(40,140){\circle*{3}}
\put(20,140){\circle*{3}}\put(29,147){$\mbox{}_{1-2m}$}
% \qbezier(250,140)(260,160)(270,140)
\put(270,10){\circle*{3}} \put(274,10){$p_0$}
\put(270,10){\line(0,1){130}} \qbezier(10,140)(80,10)(150,140)
\qbezier(150,140)(250,270)(270,140) \put(80,140){\circle{10}}
\put(68,131){$\mbox{}_{i-2m+1}$} \qbezier(28,140)(120,40)(130,140)
\qbezier(10,140)(10,160)(28,140)
\qbezier(130,140)(130,180)(80,145) \put(273,80){$u_{i}$}
\end{picture}
\begin{center}
Fig. 8
\end{center}
The element $v_i=\Delta_{m-1,m}^{-2} \Delta_{m-1,1}^{-2}(u_i)$,
$i=2,\dots, m$, is represented by the loop drawn in Fig. 9

\begin{picture}(500,200)
\put(260,140){\circle*{3}}
\put(254,131){$\mbox{}_{-1}$}
%\put(234,147){$\mbox{}_{-2}$}
%\put(240,140){\circle*{3}} \put(220,140){\circle*{3}}
%\put(200,140){\circle*{3}}
\put(205,140){$\dots$}
\put(160,140){\circle*{3}}
\put(153,131){$\mbox{}_{-m}$}
%\put(140,140){\circle*{3}}
\put(120,140){\circle*{3}} \put(103,147){$\mbox{}_{-m-1}$}
\put(95,140){$\dots$} \put(80,140){\circle*{3}}
\put(50,140){$\dots$} \put(30,140){\circle*{3}}
%$\put(20,140){\circle*{3}}
\put(24,147){$\mbox{}_{-2m}$}
% \qbezier(250,140)(260,160)(270,140)
\put(270,10){\circle*{3}} \put(274,10){$p_0$}
\put(270,10){\line(0,1){130}} \qbezier(10,140)(10,210)(80,145)
\qbezier(10,140)(10,70)(140,140) \put(80,140){\circle{10}}
\put(68,131){$\mbox{}_{i-2m+1}$}
\qbezier(140,140)(270,240)(270,140)
%\qbezier(10,140)(10,160)(28,140)
%\qbezier(130,140)(130,180)(80,145)
\put(273,80){$v_{i}$}
\end{picture}
\begin{center}
Fig. 9
\end{center}
and the element $w_i=\Delta^2_{m,0}\Delta^2_{m,m}(v_i)=b(x_i)$,
$i=2,\dots, m$, is represented by the loop drawn in Fig. 10.

\begin{picture}(500,200)
\put(50,140){$\dots$} \put(95,140){$\dots$}
\put(205,140){$\dots$}
\put(260,140){\circle*{3}} \put(254,131){$\mbox{}_{-1}$}
%\put(240,140){\circle*{3}} \put(220,140){\circle*{3}}
%\put(200,140){\circle*{3}} \put(180,140){\circle*{3}}
\put(160,140){\circle*{3}}
%\put(140,140){\circle*{3}}
\put(123,140){\circle*{3}} \put(108,147){$\mbox{}_{-m-1}$}
\put(153,131){$\mbox{}_{-m}$} \put(17,131){$\mbox{}_{-2m}$}
%\put(100,140){\circle*{3}}
\put(80,140){\circle*{3}}
%\put(60,140){\circle*{3}}
%\put(40,140){\circle*{3}}
\put(30,140){\circle*{3}} \put(270,10){\circle*{3}}
\put(274,10){$p_0$} \put(270,10){\line(0,1){130}}
\qbezier(80,145)(270,210)(270,140)
%\qbezier(10,140)(10,70)(140,140)
\put(80,140){\circle{10}}
\put(68,131){$\mbox{}_{i-2m+1}$}
%\qbezier(140,140)(270,240)(270,140)
%\qbezier(10,140)(10,160)(28,140)
%\qbezier(130,140)(130,180)(80,145)
\put(273,80){$w_{i}$}
\end{picture}
\begin{center}
Fig. 10
\end{center}
Therefore
\begin{equation} \label{second}
b(x_i)=\Delta_{2m}(x_i)
\end{equation}
for $i=2,\dots,m$.

It is easy to see that $\bar c_{m,1}(x_{2m-i})=x_{2m-i}$ for
$i=1,\dots, m-1$. The element
$y_{2m-i}=\Delta_{2(m-1),1}(x_{2m-i})$ is represented by the loop
drawn in Fig. 11.

\begin{picture}(500,200)
\put(103,147){$\mbox{}_{-m-1}$} \put(153,131){$\mbox{}_{-m}$}
\put(42,131){$\mbox{}_{-2m}$} \put(260,140){\circle*{3}}
\put(254,131){$\mbox{}_{-2}$} \put(280,140){\circle*{3}}
\put(274,131){$\mbox{}_{-1}$}

%\put(234,147){$\mbox{}_{-2}$}
%\put(240,140){\circle*{3}} \put(220,140){\circle*{3}}
\put(235,140){$\dots$} \put(200,140){\circle*{3}}
\put(175,140){$\dots$}
\put(160,140){\circle*{3}} %\put(153,131){$\mbox{}_{-m}$}
%\put(140,140){\circle*{3}}
\put(120,140){\circle*{3}} %\put(103,147){$\mbox{}_{-m-1}$}
%\put(100,140){\circle*{3}} \put(80,140){\circle*{3}}
%\put(60,140){\circle*{3}} \put(40,140){\circle*{3}}
\put(80,140){$\dots$} \put(55,140){\circle*{3}}
%\put(14,147){$\mbox{}_{-2m}$}
% \qbezier(250,140)(260,160)(270,140)
\put(290,10){\circle*{3}} \put(294,10){$p_0$}
%\put(270,10){\line(0,1){130}}
\qbezier(200,145)(250,210)(270,140)
\qbezier(290,10)(280,75)(270,140) \put(200,140){\circle{10}}
\put(195,131){$\mbox{}_{-i-1}$} \put(292,80){$y_{2m-i}$}
\end{picture}
\begin{center}
Fig. 11
\end{center}
The element $u_{2m-i}=\Delta_{m-1,m}^{-2}
\Delta_{m-1,1}^{-2}(y_{2m-i})$ is represented by the loop drawn in
Fig. 12

\begin{picture}(500,220)
\put(235,140){$\dots$}  \put(175,140){$\dots$}
\put(42,131){$\mbox{}_{-2m}$} \put(260,140){\circle*{3}}
\put(254,131){$\mbox{}_{-1}$}
%\put(234,147){$\mbox{}_{-2}$}
%\put(240,140){\circle*{3}} \put(220,140){\circle*{3}}
\put(200,140){\circle*{3}}
%\put(180,140){\circle*{3}}
\put(160,140){\circle*{3}} \put(153,131){$\mbox{}_{-m}$}
%\put(140,140){\circle*{3}}
\put(120,140){\circle*{3}} \put(103,147){$\mbox{}_{-m-1}$}
\put(80,140){$\dots$}
%\put(100,140){\circle*{3}}
%\put(80,140){\circle*{3}} \put(60,140){\circle*{3}} \put(40,140){\circle*{3}}
\put(55,140){\circle*{3}}
%\put(14,147){$\mbox{}_{-2m}$}
% \qbezier(250,140)(260,160)(270,140)
\put(270,10){\circle*{3}} \put(274,10){$p_0$}
%\put(270,10){\line(0,1){130}}
\qbezier(200,145)(100,240)(140,140)
\qbezier(140,140)(175,75)(270,10)
%\qbezier(10,140)(10,70)(140,140)
\put(200,140){\circle{10}} \put(195,131){$\mbox{}_{-i-1}$}
%\qbezier(140,140)(270,240)(270,140)
%\qbezier(10,140)(10,160)(28,140)
%\qbezier(130,140)(130,180)(80,145)
\put(193,80){$u_{2m-i}$}
\end{picture}
\begin{center}
Fig. 12
\end{center}
and the element
$v_{2m-i}=\Delta^2_{m,0}\Delta^2_{m,m}(u_{2m-i})=b(x_{2m-i})$
is represented by the loop drawn
in Fig. 13.

\begin{picture}(500,200)
\put(235,140){$\dots$}  \put(175,140){$\dots$}
\put(75,140){$\dots$} \put(42,131){$\mbox{}_{-2m}$}
\put(260,140){\circle*{3}} \put(254,131){$\mbox{}_{-1}$}
%\put(234,147){$\mbox{}_{-2}$}
%\put(240,140){\circle*{3}} \put(220,140){\circle*{3}}
\put(200,140){\circle*{3}}
%\put(180,140){\circle*{3}}
\put(160,140){\circle*{3}} \put(153,131){$\mbox{}_{-m}$}
%\put(140,140){\circle*{3}}
\put(120,140){\circle*{3}} \put(103,147){$\mbox{}_{-m-1}$}
%\put(100,140){\circle*{3}} \put(80,140){\circle*{3}}
%\put(60,140){\circle*{3}} \put(40,140){\circle*{3}}
\put(55,140){\circle*{3}}
%\put(14,147){$\mbox{}_{-2m}$}
% \qbezier(250,140)(260,160)(270,140)
\put(270,10){\circle*{3}} \put(274,10){$p_0$}
\put(270,10){\line(0,1){130}} \qbezier(200,145)(270,210)(270,140)
%\qbezier(10,140)(10,70)(140,140)
\put(200,140){\circle{10}} \put(195,131){$\mbox{}_{-i-1}$}
%\qbezier(140,140)(270,240)(270,140)
%\qbezier(10,140)(10,160)(28,140)
%\qbezier(130,140)(130,180)(80,145)
\put(273,80){$y_{2m-i}$}
\end{picture}
\begin{center}
Fig. 13
\end{center}

Therefore, by (\ref{delta}),
\begin{equation} \label{third}
b(x_{2m-i})=\Delta_{2m}(x_{2m-i})
\end{equation}
for $i=1,\dots,m-1$.

The element $y_{2m}=\bar c_{m,1}(x_{2m})$ is represented by the
loop drawn in Fig. 14.

\begin{picture}(500,220)

\put(260,140){\circle*{3}}
\put(200,140){$\dots$}
\put(160,140){\circle*{3}}
\put(153,147){$\mbox{}_{-m}$}
\put(120,140){\circle*{3}}
\put(103,131){$\mbox{}_{-m-1}$}
%\put(100,140){\circle*{3}} \put(80,140){\circle*{3}}
\put(75,140){$\dots$} \put(40,140){\circle*{3}}
\put(20,140){\circle*{3}} \put(31,147){$\mbox{}_{1-2m}$}
\put(13,131){$\mbox{}_{-2m}$}
% \qbezier(250,140)(260,160)(270,140)
\put(270,10){\circle*{3}} \put(274,10){$p_0$}
%\put(270,10){\line(0,1){130}}
\qbezier(260,135)(260,90)(140,140) \qbezier(30,140)(70,75)(270,10)
\qbezier(30,140)(-10,200)(140,140)
%\qbezier(10,140)(10,70)(140,140)
\put(260,140){\circle{10}} \put(255,147){$\mbox{}_{-1}$}
\put(120,80){$y_{2m}$}
\end{picture}
\begin{center}
Fig. 14
\end{center}
The element $u_{2m}=\Delta_{2(m-1),1}(y_{2m})$ is represented by
the loop drawn in Fig. 15.

\begin{picture}(500,200)
\put(48,131){$\mbox{}_{-2m}$} \put(200,140){$\dots$}
\put(85,140){$\dots$} \put(260,140){\circle*{3}}
\put(254,131){$\mbox{}_{-1}$}
%\put(234,147){$\mbox{}_{-2}$}
%\put(240,140){\circle*{3}} \put(220,140){\circle*{3}}
%\put(200,140){\circle*{3}} \put(180,140){\circle*{3}}
\put(160,140){\circle*{3}} \put(153,131){$\mbox{}_{-m}$}
%\put(140,140){\circle*{3}}
\put(120,140){\circle*{3}} \put(103,147){$\mbox{}_{-m-1}$}
%\put(100,140){\circle*{3}} \put(80,140){\circle*{3}}
%\put(60,140){\circle*{3}} \put(40,140){\circle*{3}}
\put(55,140){\circle*{3}}
%\put(14,147){$\mbox{}_{-2m}$}
% \qbezier(250,140)(260,160)(270,140)
\put(270,10){\circle*{3}} \put(274,10){$p_0$}
%\put(270,10){\line(0,1){130}}
\qbezier(260,145)(100,240)(140,140)
\qbezier(140,140)(175,75)(270,10)
%\qbezier(10,140)(10,70)(140,140)
\put(260,140){\circle{10}} %\put(155,131){$\mbox{}_{-1}$}
%\qbezier(140,140)(270,240)(270,140)
%\qbezier(10,140)(10,160)(28,140)
%\qbezier(130,140)(130,180)(80,145)
\put(193,80){$u_{2m-i}$}
\end{picture}
\begin{center}
Fig. 15
\end{center}
It is easy to see that $\Delta_{m-1,m}^{-2}
\Delta_{m-1,1}^{-2}(u_{2m})=u_{2m}$  and the element
$b(x_{2m})=\Delta^2_{m,0}\Delta^2_{m,m}(u_{2m})=x_1$. Therefore,
by (\ref{delta}),
\begin{equation} \label{four}
b(x_{2m})=\Delta_{2m}(x_{2m}).
\end{equation}
and Lemma \ref{Delta} follows from (\ref{first}) -- (\ref{four}).
\qed \vspace{0.1cm}
\newline{\bf Full twist formula for double number of strings.}
\begin{equation} \label{doubling formula}
\Delta^2_{2m}=\Delta^4_{m,0}\Delta^4_{m,m}r_{m}^2.
\end{equation}

\proof By Lemma \ref{Delta},
$$\Delta_{2m}=\Delta^2_{m,0}\Delta^2_{m,m}\bar r_{m}.$$
It is easy to see that $[\Delta_{m,0},\Delta_{m,m}]=\mathbf 1$ and
$[\Delta_{2m},\Delta^2_{m,0}\Delta^2_{m,m}]=\mathbf 1$. Therefore
$[\bar r_{m,m},\Delta^2_{m,0}\Delta^2_{m,m}]=\mathbf 1$ and
$$\Delta^2_{2m}=\Delta^4_{m,0}\Delta^4_{m,m}\bar r_{m}^2.$$
Thus
$$\Delta^2_{2m}=\Delta_{m,m}^{-1}\Delta^2_{2m}\Delta_{m,m}
=\Delta^4_{m,0}\Delta^4_{m,m}r_{m}^2.$$ \qed

\section{Factorization  semigroups over braid groups \label{semi}}
Let $\{ g_i \}_{i\in I}$ be a set of elements of the braid group
$Br_m$. For each $i\in I$ denote by $O_{g_i}\subset Br_m$ the set
of the elements in $Br_m$ conjugated to $g_i$ (the orbit of $g_i$
under the action of $Br_m$ by the inner automorphisms). The union
$U= \cup_{i\in I} O_{g_i}\subset Br_m$ is called the {\it full set
of conjugates} of $\{ g_i\}_{i\in I}$ and the pair $(Br_m,U)$ an
{\it equipped braid group}. In \cite {Kh-Ku}, a semigroup
$S(Br_m,U)$, called  a  {\it factorization semigroup over $Br_m$},
was associated to each equipped braid group $(Br_m,U)$. Recall
that by definition, the semigroup $S(Br_m,U)$ is a semigroup
generated by the alphabet $U$ being subject to the relations
$$ u_1\cdot u_2=u_2\cdot(u_2^{-1}u_1u_2)\quad\mbox{if $u_2\ne\mathbf 1$ and} \,\,
u_1\cdot\mathbf 1=u_1\quad\mbox{otherwise};
$$
$$ u_1\cdot u_2=(u_1u_2u_1^{-1})\cdot\,u_1 \quad\mbox{if $u_1\ne\mathbf 1$ and} \,\,
\mathbf 1\cdot u_2=u_2\quad\mbox{otherwise}
$$
for all $u_1, u_2\in U$.

There are two natural homomorphisms: the {\it product
homomorphism} $\alpha =\alpha_U :S(Br_m,U)\to Br_m$ given by
$\alpha(u)=u$ for each $u\in U$ and the homomorphism $\lambda :
Br_m\to \mbox{Aut} (S(Br_m,U))$ (the {\it  conjugation action})
given by $\lambda (g) (u)=gug^{-1}\in U$ for all $g\in Br_m$. The
action $ \lambda (g)$ on $S(Br_m,U)$ is called the {\it simultaneous
conjugation } by $g$. Put $\rho (g)=\lambda(g^{-1})$, $\lambda
_{S} =\lambda \circ\alpha_U $ and $\rho _{S} =\rho \circ\alpha_U$.
The orbit of an element $s\in S(Br_m,U)$ under the conjugation
action of $Br_m$ on $S(Br_m,U)$ is called the {\it type} of $s$.

Notice that $S: (Br_m,U)\mapsto (S(Br_m,U), \alpha_U, \lambda)$ is
a functor from the category of equipped braid groups to the
category of the semigroups over braid groups. In particular, if
$U\subset V$ are two full sets of conjugates in $Br_m$, then the
identity map $id: Br_m\to Br_m$ defines an embedding
$id_{U,V}:S(Br_m,U)\to S(Br_m,V)$. So that, for each group $Br_m$,
the semigroup $S_{Br_m}=S(Br_m,Br_m)$ is the {\it universal
factorization semigroup over} $Br_m$, which means that each
factorization semigroup $S_U=S(Br_m,U)$ over $Br_m$ is canonically
embedded in $S_{Br_m}$ by $id_{U,Br_m}$.

Since $\alpha_U=\alpha_{Br_m} \circ id_{U,Br_m}$, there is no
difference between the product homomorphisms $\alpha_U$ and
$\alpha_{Br_m}$, so the both are denoted simply by $\alpha$.

For any $s_1,\, s_2\in S(Br_m,U)$ we have
\begin{equation}
\label{change} s_1\cdot s_2=s_2\cdot \rho_S(s_2)(s_1)
=\lambda_S(s_1)(s_2)\cdot s_1.
\end{equation}

We say that an element $b\in Br_m$ has {\it the interlacing
number} $l(b)=k$ if $k$ is the smallest number such that $b$ is
conjugated in $B_m$ to an element $\bar b\in B_{k,0}$, $b= g\bar
bg^{-1}$. The element $\bar b$ will be called a {\it standard
form} of $b$.

For each $s\in S_U$ denote by $B_s$ the subgroup of $Br_m$
generated by the images $\alpha (u_1),\dots ,\alpha (u_n)$ of the
elements $u_1, \dots , u_n$ of a factorization $s=u_1\cdot ...
\cdot u_n$. It is easy to see that the subgroup $B_s$ of $Br_m$
does not depend on the presentation of $s\in S_U$ as a word in
letters $u_i\in U$.

Associate to $s\in S_{B_m}$ one more group
\begin{equation} \label{conj}
\begin{array}{l}
G(s)=G(W_s)=G(W(B_s,\mathbb F_m))\simeq \\
<x_1, \dots , x_m \, \mid \, x_i^{-1}b(x_i)=\mathbf 1, \, \,
i=1,\dots, m,\, \, b\in B_s\,
>.
\end{array}
\end{equation}

\begin{claim} \label{add} For $s_1,s_2\in S_{Br_m}$
\begin{itemize}
\item[(i)] there is a natural $C$-epimorphism $$\psi_{s_1} :G(s_2)\to
G(s_1\cdot s_2)\simeq
G(s_2)/N_{G(s_2)}(\varphi_{W_{s_2}}(W_{s_1});$$
\item[(ii)] if $s_2=\lambda (b)(s_1)$, then there is a $C$-isomorphism
$\gamma_b:G(s_1)\to G(s_2)$, in particular, if $s_1=s_2$ then
$G(s_1)=G(s_2)$ and $\gamma_{\mathbf 1}=\mbox{Id}$.
\end{itemize}
\end{claim}

\proof Straightforward.\qed

Let a factor $u_i$ of $s=u_1\cdot ... \cdot u_n$ have the
interlacing number $l(u_i)=k_i$ and a standard form $\bar u_i\in
Br_{k,0}$, $u_i=g_i\bar u_ig_i^{-1}$. Denote by $G_{u_i,loc}$ a
subgroup of $G(s)$ generated by $\varphi_{W_s}(g_i(x_j))$,
$j=1,\dots,k_i$. It follows from (\ref{change}) that for each
factor $u_i$ of $s=u_1\cdot ... \cdot u_n$ the subgroup
$G_{u_i,loc}$ is defined uniquely up to conjugation in $G(s)$.

The embedding $Br_m=B_{m,0}\subset Br_{2m}$ induces an embedding
$S_{Br_m}\subset S_{Br_{2m}}$. Consider the homomorphism
$\mbox{sh}=\mbox{sh}_m:Br_m\to Br_{2m}$ given by
$\mbox{sh}(a_i)=a_{m+i}$ for $i=1,\dots,m-1$. It induces an
embedding $\mbox{sh}:S_{Br_m}\to S_{Br_{2m}}$. Put also
$$\widetilde
r_{m}=c_{m,m}\cdot ... \cdot c_{m,1}\in S_{Br_{2m}},$$ where
$c_{m,i}$ were defined in (\ref{ci}). Applying  formula
(\ref{doubling formula}), we obtain \vspace{0.1cm}
\newline {\bf Full twist factorization formula for double number of strings.}
{\it For four elements $s_1,\dots, s_4\in S_{Br_m}$ such
that $\alpha (s_i)=\Delta_m^2$, the element
\begin{equation} \label{doubling formula 4s}
\overline s=d(s_1,s_2,s_3,s_4)=s_1\cdot s_2\cdot
\mbox{sh}(s_3)\cdot \mbox{sh}(s_4)\cdot \widetilde r_{m}\cdot
\widetilde r_{m}
\end{equation}
is a factorization of $\Delta^2_{2m}$ in $S_{Br_{2m}}$, i.e.,
$\alpha (\overline s)=\Delta^2_{2m}$.}\vspace{0.1cm}

If $s_1=s_2=s_3=s_4=s$ the doubling $d(s,s,s,s)$ will be denoted
simply by
\begin{equation} \label{doubling formula s} d(s)=d(s,s,s,s)=s\cdot s\cdot
\mbox{sh}(s)\cdot \mbox{sh}(s)\cdot \widetilde r_{m}\cdot
\widetilde r_{m}.
\end{equation}
and we put
\begin{equation}
d^{n+1}(s)=d(d^n(s)).
\end{equation}

\begin{claim} \label{cl3}
For $s\in S_{Br_m}$ such that $\alpha (s)=\Delta_m^2$, the groups
$G(s)$ and $G(d(s))$ are $C$-isomorphic.
\end{claim}
\proof It follows from Claim \ref{dclaim2}.\qed

\begin{lem} \label{new}
Let $s\in S_{Br_M}$ be such that $\alpha (s)=\Delta_M^2$. Consider
an element $g\in Br_M$ such that $g(x_2)=x_i$ for some $i\leq M$.
Put $y=g(x_1)\in \mathbb F_M$ and $b=ga_1g^{-1}$.  Then
\begin{itemize}
\item[(i)] the $C$-group $\overline G=G(s)/N_{G(s)}(\varphi_{W_s}(x_i^{-1}y))$ is
canonically $C$-isomorphic to $G(\overline s)$, where
$$\overline s=d(\lambda_S(b)(d(s)),d(s),d(s),d(s)).$$
\item[(ii)] the $C$-group $\widetilde G=G(s)/N_{G(s)}(\varphi_{W_s}([x_iy]))$ is
canonically $C$-isomorphic to $G(\widetilde s)$, where
$$\widetilde s=d(\lambda_S(b^2)(d(s)),d(s),d(s),d(s)).$$
\end{itemize}
If $M=mk$ and the set of relations $R_{W_{s}}$ implies the
relations $R_{mod\, m}$ in $\mathbb F_M$, then the set of
relations $R_{W_{\overline s}}$ (respectively, $R_{W_{\widetilde
s}}$) implies the set of relations $R_{mod\, m}$ in $\mathbb
F_{4M}$.
\end{lem}
\proof By Claim \ref{cl3}, the $C$-groups $G(s)$ and $G(d^2(s))$
are canonically $C$-isomorphic. By Claim \ref{Ak}, the set of
relations $R_{W(<b^k>,\mathbb F_M)}$ is equivalent to the single
relation $\{x_i^{-1}y=\mathbf 1\}$ if $k=1$ and $\{
[x_i,y]=\mathbf 1\}$ if $k=2$. By Claim \ref{add}, the $C$-group
$\overline G$ is canonically isomorphic to $G(b\cdot d^2(s))$ and
the $C$-group $\widetilde G$ is canonically isomorphic to
$G(b^2\cdot d^2(s))$.

Put
$$s^{\prime}=d(s)\cdot \mbox{sh}_{2M}(d(s)) \cdot
\mbox{sh}_{2M}(d(s))\cdot \widetilde r_{2M}\cdot \widetilde
r_{2M}.$$
The element
$$\begin{array}{l}
b^k\cdot d^2(s)= b^k\cdot d(s)\cdot s^{\prime}=
\lambda_S(b^k)(d(s))\cdot b^k\cdot s^{\prime}= \\
\lambda_S(\lambda_S(b^k)(d(s)))(b^k)\cdot
\lambda_S(b^k)(d(s))\cdot s^{\prime}= \\
 b^k\cdot \lambda_S(b^k)(d(s))\cdot s^{\prime}=
 \left\{ \begin{array}{l}
 b\cdot \overline s\, \, \mbox{if}\, \, k=1, \\
b^2\cdot \widetilde s\, \, \mbox{if}\, \, k=2,
\end{array}
\right.
\end{array}
$$
since
$$\alpha(d(s))=\alpha(\lambda_S(b^k)(d(s)))=\Delta^2_{2M}$$
and the elements $b^k$ and $\Delta^2_{2M}$ commute in $Br_{2M}$.
Therefore, by Claim \ref{add}, there are the canonical
epimorphisms $\psi_{b}: G(\overline s)\to \overline G$ and
$\psi_{b^2}: G(\widetilde s)\to \widetilde G$. On the other hand,
the factors $d(s)\cdot \mbox{sh}_{2M}(d(s)) \cdot
\mbox{sh}_{2M}(d(s))\cdot \widetilde r_{2M}\cdot \widetilde
r_{2M}$ of the element $\overline s$ give the same set of
relations as the element $d^2(s)$ gives. In particular, these
relations imply the relations $\{ x_{M+i}=x_i\}\in R_{mod\, M}$ in
$\mathbb F_{4M}$ and, moreover, if the set of relations
$R_{W_{s}}$ implies the relations $R_{mod\, m}$ in $\mathbb F_M$,
then the set of relations $R_{W_{\overline s}}$ (respectively,
$R_{W_{\widetilde s}}$) implies the set of relations $R_{mod\, m}$
in $\mathbb F_{4M}$. The factor
$\lambda_S(b^k)(c_{M,i})=b^kc_{M,i}b^{-k}$ of
$\lambda_S(b^k)(d(s))$ gives the relation
$x_{M+i}=bc_{M,i}b^{-1}(x_{M+i})$ in $G(\overline s)$ if $k=1$ and
$x_{M+i}=b^2c_{M,i}b^{-2}(x_{M+i})$ in $G(\widetilde s)$ if $k=2$.
But
$$b^kc_{M,i}b^{-k}(x_{M+i})=b^kc_{M,i}(x_{M+i})=b^k(x_i).$$
For $k=1,2$, we have $b(x_i)=y$ and $b^2(x_i)=yx_iy^{-1}$.  The
relations $\{ x_{M+i}=y\}$ (respectively, $\{
x_{M+i}=yx_iy^{-1}\}$) and $\{ x_{M+i}=x_i\}$ imply the relation
$x_i=y$ in $G(\overline s)$ (respectively, $\{ [x_{i},y]=\mathbf
1\}$ in $G(\widetilde s)$). Therefore the canonical epimorphisms
$\psi_b: G(\overline s)\to \overline G$ and $\psi_{b^2}:
G(\widetilde s)\to \widetilde G$ are $C$-isomorphisms.
%Lemma \ref{new} is proved.
\qed

Consider subclasses $\mathcal C_{\Delta^2,m}$, $m\in \mathbb N$,
of the class $\mathcal C$ consisting of $C$-groups which are
$C$-isomorphic to $C$-groups possessing $C$-presentations of the
form
$$ G=<x_1, \dots , x_m \, \mid
\, x_i^{-1}w_{i,j,k}(\bar x)^{-1}x_jw_{i,j,k}(\bar x)=\mathbf
1,\, \, \, w_{i,j,k}(\bar x)\in W >$$ for some finite set
$W\subset\mathbb F_m$  such that the words $[x_i^{-1},x_1\dots
x_m]$ belong to $W$ for all $i=1,\dots ,m$. Put
\[ \displaystyle
\mathcal C_{\Delta^2}=\bigcup_{m\in \mathbb N}\mathcal
C_{\Delta^2,m}.\]

\begin{claim} \label{cl2} Let $s\in S_{Br_m}$ be such that $\alpha (s)=\Delta^2_m$,
then the $C$-group $G(s)\in \mathcal C_{\Delta^2,m}$. In
particular, for $s=\Delta^2_m$ the group $G(\Delta^2_m)$ possesses
the $C$-presentation
$$G(\Delta^2_m)=<\, x_1,\dots,x_m\, \mid\,
[x_i^{-1},x_1\dots x_m]=\mathbf 1\,\, \mbox{for}\, \,
i=1,\dots, m\,>.$$
\end{claim}
\proof  We have $\Delta^2_m(x_i)=(x_1\dots x_m)x_i(x_1\dots
x_m)^{-1}$. Therefore, the relations $[x_i^{-1},x_1\dots
x_m]=\mathbf 1$, $i=1,\dots, m$, belong to the set of relations in
presentation (\ref{conj}), since $\alpha (s)=\Delta^2_m$.\qed

Denote by $D_{k,m}$, $k\leq m$, the full set of conjugates of the
elements $a_1$ and $\Delta_k^2$ in $Br_m$. Put $\mathcal D_{k,m}=
S(Br_m,D_{k,m})$.

\begin{thm} \label{dgroup}
Let a $C$-group $G\in \mathcal C_{\Delta^2,m}$.
Then there is $M\in \mathbb N$ and an element $s\in \mathcal
D_{m,M}$ such that
\begin{itemize}
\item[(i)] $\alpha (s)=\Delta_M^2$;
\item[(ii)] $G$ and $G(s)$ are $C$-isomorphic.
\end{itemize}
\end{thm}

\proof Let
$$ G=<x_1, \dots , x_m \, \mid
\, x_i^{-1}w_{i,j,k}(\bar x)x_jw_{i,j,k}(\bar x)^{-1}=\mathbf
1,\, \, \, w_{i,j,k}(\bar x)\in W >,$$ for some finite set
$W\subset \mathbb F_m$ such that the words $[x_i^{-1},x_1\dots
x_m]$ belong to $W$ for all $i=1,\dots ,m$.

Numerate the words $w\in W$ so that the word $[x_i^{-1},x_1\dots
x_m]$ has the number $i$. Denote by $G_n =\mathbb F_m/N_n$, where
$N_n$ is the normal closure of the set of words $\overline W_n=\{
w_1,\dots,w_n\}$ in $\mathbb F_m$. For each $n\geq m$ we construct
an element $s_n\in \mathcal D_{m,M_n}$ such that
\begin{itemize}
\item[(1)] $\alpha (s_n)=\Delta_{M_n}^2$;
\item[(2)] $G_n$ and $G(s_n)$ are $C$-isomorphic.
\end{itemize}

We start with $M_m=m$ and $s_m=\Delta_m^2\in \mathcal D_{m,m}$. By
Claim \ref{cl2}, $G_m$ and $G(s_m)$ are $C$-isomorphic.

Assume that for some $k\geq 0$ and $p\geq 0$ we have constructed
an element $s_{m+k}\in \mathcal D_{m,M_{m+k}}$, where
$M_{m+k}=2^pm$, such that
\begin{itemize}
\item[(1)] $\alpha (s_{m+k})=\Delta_{M_{m+k}}^2$,
\item[(2)] $G_{m+k}$ and $G(s_{m+k})$ are $C$-isomorphic,
\item[(3)] the set of relations $R_{W_{s_{m+k}}}$ implies the relations
$R_{mod\, m}$.
\end{itemize}
Construct an element $s_{m+k+1}$ having the similar properties as
$s_{m+k}$ has. Put $n=m+k$ and consider the word
$$w_{n+1}=
x_{i_n}^{-1}w_{i_{n},j_n,k_n}(\bar
x)x_{j_n}w_{i_n,j_n,k_n}(\bar x)^{-1}.$$
Let $l$ be the letter
length of the word $w_{i_n,j_n,k_n}$. Consider the element
$\widetilde s_n=d^{l+2}(s_n)$.  By Claim \ref{dclaim2}, the sets
of relations $R_{W_{s_n}}$ and $R_{W_{d^{l+2}(s_n)}}$ are
equivalent. Note that $\alpha
(d^{l+2}(s_{n}))=\Delta_{2^{l+2}M_{n}}^2$ and condition ($3$) is
also realized  for the element $\widetilde s_n=d^{l+2}(s_n)$. Put
$$y_n=\mbox{m}_{2^lM_n,m}(w_{i_{n},j_n,k_n})x_{2^lM_n+j_n}
\mbox{m}_{2^lM_n,m}(w_{i_n,j_n,k_n})^{-1}.$$ By Claims \ref{move}
and \ref{pair}, there is an element $b\in Br_{2^{l+2}M_{n}}$
conjugated to the standard generator $a_1$ and such that the set
of relations $R(W_{<b>,\mathbb F_{2^{l+2}M_{n}}})$ is equivalent
to the single relation $x_{i_n+2^{l+1}M_{n}}=y_n$. Applying Lemma
\ref{new} (i), the element
$$s_{n+1}=d(\lambda_S(b)(d(\widetilde s_n)),d(\widetilde s_n),d(\widetilde s_n),
d(\widetilde s_n))$$ is a desired one which completes the proof of
the Theorem.\qed

\section{Weak $\mu$-equivalence}

Denote by $A_i=A_i(m)$ the full set of conjugates of the element
$a_1^{i+1}$ in $Br_m$, where $\{ a_1,\dots, a_{m-1}\}$ is a set of
standard generators of $Br_m$ (recall that all generators
$a_1,\dots ,a_{m-1}$ are conjugated to each other in $Br_m$), and
put $\mathcal A_k= S(Br_m,A_{\leq k})$, $\mathcal A^0_k=
S(Br_m,A^0_{\leq k})$, where $A_{\leq k}=A_{-3}\cup
(\cup^k_{i=-1}A_i)$ and $A^0_{\leq k}=\cup^k_{i=0}A_i$. We have
the embedding $\mathcal A_k^0\subset \mathcal A_k$.

Below, we restrict our consideration to the case $k=2$. The
semigroup $\mathcal A^0_2$ is called the {\it braid monodromy
cuspidal factorization semigroup} and $\mathcal A_2$  the {\it
braid monodromy cuspidal factorization semigroup with negative
nodes}.

Let $U\subset Br_m$ be a full set of conjugates containing
$A_2\cup A_{-3}$. Note that if $g\in A_1$ then $g^{-1}\in A_{-3}$.
Consider a semigroup
$$\overline{S}_U = \langle\, u\in U)\,\,
\mid \,\, R\in \mathcal{R}\cup \overline{\mathcal R} \, \rangle,
$$ where $\mathcal R$ is the set of relations defining $S_U$ (see
section \ref{semi}) and
$$\overline{\mathcal R}=\{ \, g\cdot (g^{-1})=(g^{-1})\cdot g=
\mathbf 1\, \mid \, g\in A_1\, \}$$ is the set of {\it
cancellation relations}. There is the canonical homomorphism of
semigroups
$$c:\mathcal S_U\to \overline{S}_U.$$
We say that two elements $s_1,\, s_2\in S_U$ are {\it weakly
equivalent} if $c(s_1)=c(s_2)$. It is easy to see that
$\overline{S}_U$ can be considered as a semigroup over $Br_m$ with
the product homomorphism $\overline \alpha :\overline{S}_U\to
Br_m$ such that $\alpha =\overline \alpha \circ c$.

Let $(s,\mu_s)$ be a pair consisting of $s\in S_U$ and a
homomorphism $\mu_s: G(s)\to \mathfrak S_N$ to the symmetric group
$\mathfrak S_N$ for some $N$. We say that two pairs
$(s_1,\mu_{s_1})$ and $(s_2,\mu_{s_2})$ are {\it equivalent} if
$s_1$ and $s_2$ belong to the same orbit under the conjugation
action of $Br_m$ on $S_U$ (i.e., there is $b\in Br_m$ such that
$s_2=\lambda (b)(s_1)$), and $\mu_2=\mu_1\circ \gamma_b$, where
$\gamma_b$ is defined in Claim \ref{add}.

Two pairs $(s_1,\mu_{s_1})$ and $(s_2,\mu_{s_2})$ is said to be
obtained by {\it admissible transformation} from each other
(notation: $(s_1,\mu_{s_1})\leftrightarrow (s_2,\mu_{s_2})$) if
there is $g\in A_1$ such that $s_2=g\cdot (g^{-1})\cdot s_1$,
$\mu_1=\psi_{g\cdot(g^{-1})}\circ \mu_2$, and the order of the
group $\mu_{s_2}(G_{g,loc})$ is greater than two. We say that two
pairs $(s^{\prime},\mu_{s^{\prime}})$ and
$(s^{\prime\prime},\mu_{s^{\prime\prime}})$ are {\it weakly
$\mu$-equivalent} if  there is a sequence of pairs
$(s_1,\mu_1),\dots ,(s_n,\mu_n)$ such that
$(s^{\prime},\mu_{s^{\prime}})=(s_1,\mu_1)$,
$(s^{\prime\prime},\mu_{s^{\prime\prime}})=(s_n,\mu_n)$, and for
each $i=1,\dots, n-1$ the pairs $(s_i,\mu_i)$ and
$(s_{i+1},\mu_{i+1})$ are either equivalent or can be obtained
from each other by an admissible transformation.

Below, we consider the case then $U=\mathcal A_2$. Let $s=u_1\cdot
...\cdot u_n\in \mathcal A_2$, $u_i\in A_{\leq 2}$. We say that a
homomorphism $\mu_s:G(s)\to \mathfrak S_N$ is {\it generic} if
\begin{itemize}
\item[(i)] $\mu_s$ is an epimorphism,
\item[(ii)] $\mu_s(\varphi_s(x_i))$ is a transposition in
$\mathfrak S_N$ for each $C$-generator $\varphi_s(x_i)$ of the
group $G(s)$,
\item[(iii)] the order of the group $\mu_s(G_{u_i,loc})$ is
greater than two for each factor $u_i\in A_{-3}\cup A_1\cup A_2$.
\end{itemize}

Let $(s,\mu_s)$ be a pair consisting of $s\in \mathcal A_2$ and a
homomorphism $\mu_s: G(s)\to \mathfrak S_N$ for some $N$. By Claim
\ref{cl3}, there is the canonical $C$-isomorphism $\psi _d:
G(d(s))\to G(s)$. It defines the pair $(d(s),d(\mu_s))$, where
$d(\mu_s)=\mu_s\circ \psi_d$.

\begin{claim}\label{dmu}
Let $s\in \mathcal A_2$ and $\mu_s:G(s)\to \mathfrak S_N$ be a
generic epimorphism. Then the homomorphism $d(\mu_{s}):G(d(s))\to
\mathfrak S_N$ is generic.
\end{claim}

\proof Straightforward.\qed

\begin{thm} \label{alg} Let $s=u_1\cdot ...\cdot u_n\in \mathcal
A_2\subset S_{Br_m}$, $u_i\in A_{\leq 2}$, such that $\alpha
(s)=\Delta^2_m$, and $\mu_s:G(s)\to \mathfrak S_N$ be a generic
epimorphism. Let $x_i$ be a generator of $\mathbb F_m$ and $y$ be
an element conjugated in $\mathbb F_m$ to a generator $x_j$ for
some $j$ such that $\mu_s(\varphi_s(x_i))$ and
$\mu_s(\varphi_s(y))$ are two different commuting transpositions
in $\mathfrak S_N$. Then there are $M\in \mathbb N$ and two pairs
$(\overline s,\mu_{\overline s})$ and $(\widetilde
s,\mu_{\widetilde s})$, $\overline s,\, \widetilde s\in \mathcal
A_2\subset S_{Br_M}$, such that
\begin{itemize}
\item[(i)] $\alpha (\overline s)=\alpha (\widetilde s)=\Delta^2_M$,
\item[(ii)] $G_{\overline s}$ is $C$-isomorphic to $G(s)$ and
$G_{\widetilde s}$ is $C$-isomorphic to the $C$-group $G(s)/\{
[\varphi_s(x_i),\varphi_s(y)]=\mathbf 1\}$,
\item[(iii)] if $[\varphi_s(x_i),\varphi_s(y)]\neq \mathbf 1$,
then the types of $\overline s$ and $\widetilde s$ are different,
\item[(iv)] $\mu_{\overline s}$ and $\mu_{\widetilde s}$
are generic epimorphisms onto $\mathfrak S_N$,
\item[(v)] the pairs $(\overline s,\mu_{\overline s})$ and
$(\widetilde s,\mu_{\widetilde s})$ are weakly $\mu$-equivalent.
\end{itemize}
Moreover, if $s\in \mathcal A^0_2\subset S_{Br_m}$, then
$\overline s,\, \widetilde s\in \mathcal A^0_2\subset S_{Br_M}$.
\end{thm}

\proof By Claims \ref{move}, \ref{main1}, \ref{pair}, and
\ref{dmu}, there are
\begin{itemize}
\item[(1)] $M_1=2^nm$ (for
some $n\in \mathbb N$),
\item[(2)] an element $y_1\in \mathbb
F_{M_1}$ which is a move apart of the element $y$ such that
$(y_1,x_k)$ is a good pair for some $k\equiv i \mod m$,
\item[(3)] an element $s_1=d^n(s)\in \mathcal A_2\subset
S_{Br_{M_1}}$ such that $\mu_s(x_i)=\mu_{s_1}(x_k)$,
$\mu_s(y)=\mu_{s_1}(y_1)$, where $\mu_{s_1}=d^n(\mu_s)$.
\end{itemize}
By definition of the good pairs, there is an element $g\in
Br_{M_1}$ such that $g(x_1)=y_1$ and $g(x_k)=x_2$. Put
$b=ga_1g^{-1}\in Br_{M_1}$, $\overline s =d^2(s_1)$, and
$\widetilde s=d(\lambda_S(b^2)(d(s)),d(s),d(s),d(s))$. By Claim
\ref{dmu}, the groups $G(s)$ and $G(\overline s)$ are canonically
$C$-isomorphic. Therefore the generic epimorphism $\mu_s$ induces
a generic epimorphism $\mu_{\overline s}: G(\overline s)\to
\mathfrak S_N$. By Lemma \ref{new} (ii) (see its proof), the
groups $G(\widetilde s)$ and $G(b^{-2}\cdot b^2\cdot \overline s)$
are canonically $C$-isomorphic to $G(\overline s)/N_{G(\overline
s)}(\varphi_{W_{\overline s}}([x_j,y_1]))$, where $\widetilde
s=d(\lambda_S(b^2)(d(s_1)),d(s_1),d(s_1),d(s_1))$. Therefore these
isomorphisms induce generic epimorphisms $\mu_{b^{-2}\cdot
b^2\cdot \overline s}:G(b^{-2}\cdot b^2\cdot \overline s) \to
\mathfrak S_N$ and $\mu_{\widetilde s}:G(\widetilde s) \to
\mathfrak S_N$. We have $(\overline s, \mu_{\overline
s})\leftrightarrow (b^{-2}\cdot b^2\cdot \overline s,
\mu_{b^{-2}\cdot b^2\cdot \overline s})$, the pairs $(b^{-2}\cdot
b^2\cdot \overline s, \mu_{b^{-2}\cdot b^2\cdot \overline s})$ and
$(b^{-2}\cdot b^2\cdot \widetilde s, \mu_{b^{-2}\cdot b^2\cdot
\widetilde s})$ are equivalent, and $(b^{-2}\cdot b^2\cdot
\widetilde s, \mu_{b^{-2}\cdot b^2\cdot \widetilde
s})\leftrightarrow (\widetilde s, \mu_{\widetilde s})$. Therefore
the pairs $(\overline s, \mu_{\overline s})$ and $(\widetilde s,
\mu_{\widetilde s})$ are weakly $\mu$-equivalent. This completes
the proof of Theorem \ref{alg}.\qed

\section{Hurwitz curves}\label{Hurw} Let $F_N$ be a relatively minimal
ruled rational surface, $N\geq 1$, $\mbox{pr}: F_N \to \mathbb
C\mathbb P^1$ the ruling, $R$ a fiber of $\mbox{pr}$, and $E_N$
the exceptional section, $E_N^2=-N$. Below, we will identify
$\mbox{pr}: F_1 \to \mathbb C\mathbb P^1$ with a linear projection
$\mbox{pr}:\mathbb C\mathbb P^2\to \mathbb C\mathbb P^1$  with
center at a point $p\in \mathbb C\mathbb P^2$ ($p$ is the blow
down $E_1$ to the point).

\begin{df} The image $\bar H=f(S)\subset F_N$ of a smooth map
$f:S \to F_N\setminus E_N$ of an oriented closed real surface $S$
is called a Hurwitz curve {\rm (}with respect to {\rm
$\mbox{pr}$)} of degree $m$ if there is a finite subset
$Z\subset\bar H$ such that:
\begin{itemize}
\item[(i)] $f$ is an embedding  of the surface $S\setminus
f^{-1}(Z)$ and for any $s\notin Z$, $\bar H$ and the fiber
$R_{\mbox{pr}(s)}$ of $\pr$ meet at $s$ transversely and with
positive intersection number;
\item[(ii)] for each $s\in Z$ there is a neighborhood
$U\subset F_N$ of $s$ such that $\bar H\cap U$ is a complex
analytic curve, and the complex orientation of $\bar H\cap U
\setminus \{s\}$ coincides with the orientation transported from
$S$ by $f$;
\item[(iii)]  the restriction of
{\rm $\mbox{pr}$} to $\bar H$ is a finite map of degree $m$.
\end{itemize}
\end{df}

For any Hurwitz curve $\bar H$ there is one and only one minimal
$Z\subset\bar H$ satisfying the conditions from Definition 3.1. We
denote it by $Z(\bar H)$.

A Hurwitz curve $\bar H$ is called {\it cuspidal } if for each
$s\in Z(\bar H)$ there is a neighborhood $U$ of $s$ and local
analytic coordinates $z,w$ in $U$ such that
\begin{itemize}
\item[(iv)] $\pr_{\mid U}$ is given by
$(z,w)\mapsto z$;
\item[(v)]  $\bar H\cap U$ is given by
$w^2=z^k, k\ge 1$.
\end{itemize}
It is called {\it ordinary cuspidal} if $k\leq 3$ in (v) for all
$s\in Z(\bar H)$, and {\it nodal} if $k\leq 2$. Below, we will
consider Hurwitz cuspidal curves $\bar H$ {\it with negative
nodes} allowing $\bar H$ to have also singularities of the
following form:
\begin{itemize}
\item[(vi)] in a neighborhood $U\subset F_N$ of a point $s\in \bar H$ the intersection
$U\cap \bar H$ splits into two branches meeting transversally at
$s$ with negative intersection number and meeting transversally
the fibre $R_{\mbox{pr}(s)}$ with positive intersection number.
\end{itemize}

Since $E_N\cap \bar H=\emptyset$, one can define a braid monodromy
factorization $bmf(\bar H)\in S_{Br_m}$ of $\bar H$ as in the
algebraic case (see, for example, \cite{Kh-Ku}). To do this, we
fix a fiber $R_{\infty}$ meeting transversely $\bar H$ and
consider the affine Hurwitz curve $H=\bar H\cap \mathbb C^2$,
where $\mathbb C^2=F_N \setminus (R_{\infty}\cup E_N)$. Choose
$r>>1$ such that $\mbox{pr} (Z)$ belongs to the disc $D(r)=\{\,
\mid z\mid \leq r\, \} \subset\mathbb C= \mathbb C\mathbb
P^1\setminus\mbox{pr} (F_\infty )$, $Z=Z(\bar H)$. Denote by
$z_1,\dots,z_n$ the elements of $\mbox{pr} (Z)$. Pick $\varepsilon
$, $0<\varepsilon <<1$, such that the discs $D_{i}(\varepsilon
)=\{ z\in \mathbb C \, \mid \,\, \mid z-z_i\mid \leq \varepsilon
\, \}$, $i=1,\dots, n$, are disjoint. Select arbitrary points
$u_i\in
\partial D_{i}(\varepsilon )$ and a point $u_0\in \partial D(r)$. Choose
disjoint simple paths $l_i\subset D(r)\setminus \bigcup^n_1
D_{i}(r)$, $i=1,\dots,n$, starting at $u_0$ and ending at $u_i$
and renumber the points in a way that the product
$\gamma_1\dots\gamma_n$ of the loops $\gamma_i=l_i\circ
\partial D_{i}(\varepsilon )\circ l_i^{-1}$ is equal to $\partial
D(r)$ in $\pi_1(D(r)\setminus \{ z_1, \dots ,z_n\}, u_0)$.

Each $\gamma _i$, lifted to $H$, defines an element $b_i\in Br_m$.
The factorization $b_1\cdot .\, .\, . \cdot b_n\in S_{Br_m}$ is
called a {\it braid monodromy factorization} of $\bar H$. In fact,
each $b_i$ is conjugated to a braid monodromy standard form of
some algebraic germ over $z_i$ or $b_i\in A_{-3}$ if the
corresponding singularity of $\bar H$ is a negative node. Hence,
$bmf(\bar H)=b_1\cdot .\, .\, .\cdot b_n$ belongs to
$\overline{\mathcal P}=S(Br_m,P\cup A_{-3})$, where $P$ is the set
of all braid monodromies of the germs of plane algebraic curves of
degree $m$. The orbit of $bmf(\bar H)$ under the conjugation
action of $Br_m$ on $\overline{\mathcal P}$ is called the {\it
braid monodromy factorization type} of $\bar H$ and denoted by $bmt(\bar
H)$.

If $\bar H$ is an ordinary cuspidal Hurwitz curve with negative
nodes, then its braid monodromy factorization $bmf(\bar H)\in
\mathcal A_2$.

The following lemma is well known (see, for example,
\cite{Kh-Ku}).
\begin{lem} For a Hurwitz curve $\bar H\subset F_N$ of degree $m$, we have
$$\alpha (bmf(\bar H))=\Delta_m^{2N}.$$
\end{lem}

The converse statement can be also proved in a straightforward
way.

\begin{thm} \label{moisha} {\rm (\cite{Moi2})} For any $b=b_1\cdot
 .\, .\, . \cdot b_n\in \overline{\mathcal P}$ such that $\alpha (b)=\Delta_m^{2N}$
there is a Hurwitz curve $\bar H\subset F_N$ with a braid
monodromy factorization $bmf(\bar H)$ equal to $b$.
\end{thm}

\begin{df} Two Hurwitz  curves $\bar H_1$ and $\bar H_2\subset F_N$ are
called  $H$-isotopic if there is a fiberwise continuous isotopy
$\phi_t:F_N\to F_N$, $t\in [0,1]$, smooth outside the
fibers $R_{\pr(s)}$, $s\in Z(\bar H_1)$, and %
such that
\begin{itemize}
\item[(i)] $\phi_0=\text {Id}$;
\item[(ii)] $\phi_t(\bar H_1)$ is a Hurwitz  curve for all $t\in [0,1]$;
\item[(iii)] $\phi_1(\bar H_1)=\bar H_2$;
\item[(iv)]
$\phi_{t}= \mbox{Id}$ in a neighborhood of $E_N$ for all $t\in
[0,1]$.
\end{itemize}
\end{df}

By Theorem 3.2 in \cite{Kh-Ku}, two cuspidal Hurwitz curves with
negative nodes are $H$-isotopic if and only if they have the same
braid monodromy factorization type.

Obviously, if $\bar H_1$ and $\bar H_2$ are $H$-isotopic, then the
fundamental groups  $\pi_1(F_N\setminus \bar H_1)$ and
$\pi_1(F_N\setminus \bar H_2)$ are isomorphic.

Consider two sections $H_1=\{ u+1=0,\, v=0 \}$ and $H_2=\{ u=0,\,
v=0\}$ of the projection $\mbox{pr}_1 :D \to D_1$ of a bi-disc
$D=D_1\times D_2$, where $D_1=\{  \sqrt{x^2+y^2}  \leq 2\}$ and
$D_2=\{ \sqrt{u^2+v^2}\leq 2\}$. Let $h_t$ be an isotopy of $H_1$
given by equations
\[
\begin{array}{l}
u+1 +2t\rho (\sqrt{x^2+y^2})(y+x^2-1)=0, \\
v-t\lambda \rho (\sqrt{x^2+y^2}) y=0,
\end{array}
\]
where $\rho : \mathbb R \to \mathbb R$ is a monotone
$C^{\infty}$-function such that $\rho (s)=0$ if $ s\geq 2$ and
$\rho (s) =1$ if $s\leq 1$ and $\lambda$ is a small positive
constant. It is easy to check that
\begin{itemize}
\item[(1)] $h_t(H_1)$, $0\leq t\leq 1$,
and $H_2$ are symplectic curves in $D$ with respect to the
symplectic form $\omega= x\wedge y+u\wedge v$;
\item[(2)] the sections
$h_t(H_1)$ and $H_2$ have no common points if $0\leq t< 1/2$;
%and
\item[(3)]
the section $h_t(H_1)$ meets $H_2$ transversally at two points
with, respectively, positive and negative intersection numbers if
$1/2< t\leq 1$.
\end{itemize}
We say that the "regular homotopy" $h_t(H_1)\cup
H_2$ is the {\it creation of negative node}. The converse homotopy
is called the {\it cancellation of negative node}.

Two Hurwitz curves $\bar H^{\prime}$ and $\bar
H^{\prime\prime}\subset F_N$ are said to be obtained from each
other by an {\it admissible transformation} if there is a
neighborhood $U\subset F_N$ diffeomorphic to the bi-disc
$D=D_1\times D_2$ with local analytic coordinates $z=x+iy$ and
$w=u+iv$ such that $\mbox{pr}_{\mid U}=\mbox{pr}_1:D\to D_1$,
$$\bar H^{\prime}\cap U=\{ u+1=0,\, v=0
\}\cup \{ u=0,\, v=0\} ,$$
$$
%\begin{array}{l}
\bar H^{\prime\prime}\cap U= %\\
\{ u+1 +2%\rho (\sqrt{x^2+y^2})
(y+x^2-1)\rho =0,\, v-\rho
%(\sqrt{x^2+y^2})
y=0 \}\cup \{ u=0,\, v=0\} ,
%\end{array}
$$
and $\bar H^{\prime}\cap (F_N\setminus U)=\bar
H^{\prime\prime}\cap (F_N\setminus U)$.

It is easy to see that if $\bar H^{\prime}$ and $\bar
H^{\prime\prime}\subset F_N$ are obtained from each other by the
above admissible transformation, then
\begin{equation} \label{g}
bmf (\bar H^{\prime\prime})= g\cdot (g^{-1})\cdot bmf(\bar
H^{\prime})
\end{equation}
for some $g\in A_1$.

We say that $\bar H^{\prime}$ and $\bar H^{\prime\prime}\subset
F_N$ are {\it weakly equivalent} if there is a sequence of Hurwitz
curves $\bar H_i$, $i=1,\dots, n$, such that $\bar H^{\prime}=\bar
H_1$, $\bar H^{\prime\prime}=\bar H_n$, and for each $i=1,\dots,
n-1$ the Hurwitz curves $\bar H_i,\, \bar H_{i+1}$ either are
$H$-isotopic or are obtained from each other by an admissible
transformation.
\begin{claim} \label{hurweak} Two cuspidal
Hurwitz curves $\bar H^{\prime}$ and $\bar H^{\prime\prime}\subset
F_N$ (possibly with negative nodes) are weakly equivalent if and
only if their braid monodromy factorizations $bmf(\bar
H^{\prime})$ and $bmf (\bar H^{\prime\prime})$ are weakly
equivalent.
\end{claim}
\proof It follows from (\ref{g}) and Theorem 3.2 in \cite{Kh-Ku}.
\qed

Note that if two cuspidal Hurwitz curves $\bar H^{\prime}$ and
$\bar H^{\prime\prime}\subset F_N$ are weakly equivalent then it
is not necessary that their fundamental groups $\pi_1(F_N\setminus
\bar H^{\prime})$ and $\pi_1(F_N\setminus \bar H^{\prime\prime})$
are isomorphic.

\section{Fundamental groups of the complements of
Hurwitz curves}

The following theorem is well-known.
\newline
{\bf Zariski -- van Kampen Theorem}. {\it Let $\bar H\subset F_N$
be a Hurwitz curve of degree $m$ having braid monodromy
factorization $s=bmf(\bar H)$. Then }
\begin{itemize}
\item[(i)] {\it the fundamental group}
$$\pi_1(\mathbb C^2\setminus \bar H)\simeq G(B_s),$$
{\it where} $\mathbb C^2=F_N\setminus (E_N\cup R_{\infty})$.
\item[(ii)] {\it the fundamental group}
$$\pi_1(\mathbb C\mathbb P^2 \setminus \bar H)\simeq G(B_s)/
\{ \varphi_{W_s}(x_1\dots x_m)=\mathbf 1 \}.$$
\end{itemize}
\vspace{0.1cm}

In the case $N=1$, since $\alpha (s)= \Delta_m^2$ for $s=bmf(\bar
H)$ , we get that the relations
\begin{equation} \label{C0}
[x_i^{-1},x_1\dots x_m]=\mathbf 1,\,\, i=1,\dots,m
\end{equation}
belong to the set of defining relations of presentation
(\ref{conj}).

Denote by $\mathcal H$ the subclass of $\mathcal C$ consisting of
the fundamental groups of the complements of the affine plane
Hurwitz curves, i.e.,
$$\mathcal H=\{ \,  \pi_1(\mathbb C\mathbb
P^2\setminus (\bar H\cup R_{\infty}))\, \} ,$$ where $R_{\infty}$
is a generic line of the pencil defining the ruling
$\mbox{pr}:\mathbb C\mathbb P^2\to \mathbb C\mathbb P^1$
\begin{thm} \label{fun}
$\mathcal H=\mathcal C_{\Delta^2}$.
\end{thm}

\proof The inclusion $\mathcal H\subset \mathcal C_{\Delta^2}$
follows from  (\ref{conj}), (\ref{C0}), and Claim \ref{claim2}.

The reverse inclusion follows from
\begin{thm} \label{hgroup}
For any $C$-group $G\in \mathcal C_{\Delta^2,m}$ there is a
Hurwitz curve $\bar H\subset \mathbb C\mathbb P^2$ with the
singularities of the type $\{ w^m=z^m\}$ such that $\pi_1(\mathbb
C^2\setminus  H)$ is $C$-isomorphic to $G$, where $\mathbb
C^2=\mathbb C\mathbb P^2\setminus R_{\infty}$ and $H=\bar H\cap
\mathbb C^2$.
\end{thm}
\proof The braid monodromy of the singularity $\{w^m=z^m\}$ is
equal to $\Delta^2_m$. Therefore Theorem \ref{hgroup} follows from
Theorems \ref{dgroup} and \ref{moisha}.\qed \qed

Denote by $\mathcal T $ the subclass of $\mathcal K$ of the torus
knot groups, that is, the fundamental groups $\pi_1(S^3\setminus
K_{p,q})$ of the knots
$$K_{p,q}=\{ (z,w)\in S^3\, \mid \, z=e^{2p\pi it},\,
w=e^{2q\pi it},\, t\in [0,1]\}$$ ($p$, $q$ are co-prime, possibly
$p=q=1$), where
$$ S^3=\{ \, (z,w)\in \mathbb C^2\, \mid \,
|z|^2+|w|^2=2\, \} .$$
\begin{thm} \label{tor} $\mathcal H\cap \mathcal K=
\mathcal T$.
\end{thm}
\proof Let $\bar H \subset \mathbb C\mathbb P^2$ be a Hurwitz
curve of degree $m$ and $G$ the fundamental group $\pi_1(\mathbb
C^2\setminus H)$ of the complement of its affine part $H\subset
\mathbb C^2= \mathbb C\mathbb P^2 \setminus R_{\infty}$. It
follows from (\ref{C0}) that the element $\varphi_{W_s}(x_1\dots
x_m)$ belongs to the center of the group $G$. Therefore, by
Theorem 6.1 \cite{Bu-Zi}, $\mathcal H\cap \mathcal K\subset
\mathcal T$.

The knot group $G_{p,q}$ of the knot $K_{p,q}$ has Wirtinger
presentation
$$G_{p,q}=<x_1,\dots, x_p\, \mid x_i=(a_1\dots a_{p-1})^q(x_i),\, \,
i=1,\dots, p>.$$ We have
$$(a_1\dots a_{p-1})^{pq}=\Delta_p^{2q}.$$

Therefore there exists a Hurwitz curve $\widetilde H\subset F_q$
of degree $p$ with braid monodromy factorization $bmf(\widetilde H)=
s^p\in S_{Br_p}$, where $s=(a_1\dots a_{p-1})^q$ is one of the
generators of $S_{Br_p}$ (the braid monodromy of the singularity
given by the equation $w^p=z^q$). By Zariski -- van Kampen
Theorem, its fundamental group
$$\pi_1(\mathbb C^2\setminus H)\simeq G_{p,q},$$
where $H=\mathbb C^2\cap \widetilde H$ and $\mathbb
C^2=F_q\setminus R_{\infty}$, and $R_{\infty}$ is a fibre of
$\mbox{pr}$ meeting $\widetilde H$ transversally.

As is known, see for example \cite{Moi-Te},
$$\Delta_{p+1}^2=\prod_{l=p+1}^{2} \prod^{l-1}_{k=1} b_{k,l}^2=
\Delta_{p}^2\prod^{p}_{k=1} b_{k,p+1}^2, $$
 where $b_{k,l}=(a_{l-1}\dots a_{k+1})a_k(a_{l-1}\dots
a_{k+1})^{-1}$ for $k<l$ (the notation $\prod_k^n$ states for from
the left to the right product from $k$ to $n$).

Put
$$s_1=\prod^{p}_{k=1} (b_{k,p+1}^2)\in S_{Br_{p+1}},$$
where the product is taken in $S_{Br_{p+1}}$. The element
$s^p\cdot s_1^q\in S_{Br_{p+1}}$ is a braid monodromy
factorization of a Hurwitz curve in $F_q$, since $\alpha(s^p\cdot
s_1^q)=\Delta^{2q}_{p+1}$. It is easy to see that this curve
splits into two components one of which is $H$-isotopic to
$\widetilde H$ and another one is a section  $\widetilde C$
disjoint with $E_q$, $bmf(\widetilde H\cup \widetilde C)=s^p\cdot
s_1^q$. The curves $\widetilde H$ and $\widetilde C$ meet
transversally. Without loss of generality, we can assume that
$\widetilde C\subset F_q$ is an algebraic section. It follows from
the view of the braid monodromy factorization $bmf(\widetilde
H\cup \widetilde C)=s^p\cdot s_1^q$ that the fundamental group
$\pi_1(\mathbb C^2\setminus (\widetilde H\cup \widetilde C))$ is
canonically $C$-isomorphic to the direct product
$$\begin{array}{l}
G_{p,q}\times \mathbb F_1\simeq \\
<x_1,\dots,x_p\, \mid x_i=(a_1\dots a_{p-1})^q(x_i),\, \,
i=1,\dots, p>\times <x_{p+1}> \end{array}
$$
(the factor $\mathbb
F_1=<x_{p+1}>$ corresponds to the fundamental group $\pi_1(\mathbb
C^2\setminus \widetilde C)$).

Choose  non-homogeneous coordinates $(u,v)$ in $\mathbb
C^2=F_q\setminus (E_q\cup R_{\infty})$ such that $u=0$ is an
equation of $C=\widetilde C\cap (F_q\setminus (E_{q}\cup
R_{\infty}))$, and coordinates $(z,w)$ in the complement  of some
line $L=L_{\infty}\subset \mathbb C\mathbb P^2$. Let $f:\mathbb
C\mathbb P^2\to F_q$ be the cyclic (rational) covering of degree
$p$ given by $z=u^q,\, w=v$. It is branched along $\widetilde
C\cup E_q$. Put $\bar H=f^{-1}(\widetilde H)$ and $\bar
C=F^{-1}(\widetilde C)$. We have $L_{\infty}=f^{-1}(R_{\infty}\cup
E_q)$, $\bar H$ is a Hurwitz curve of degree $pq$, and $\bar C$ is
a line in $\mathbb C\mathbb P^2$ meeting $\bar H$ transversally.
It is easy to see that $f_*$ gives the isomorphism of
$\pi_1(\mathbb C^2\setminus (\bar H\cup\bar C))$ and
$G_{p,q}\times <x_{p+1}^q>$, where $<x_{p+1}^q>$ is generated by a
circuit around $\bar C$. Therefore the canonical epimorphism
$$\psi: \pi_1(\mathbb C^2\setminus (\bar H\cup\bar C))\to
\pi_1(\mathbb C^2\setminus \bar H)$$ having $\ker \psi
=<x_{p+1}^q>$ gives an isomorphism between $\pi_1(\mathbb
C^2\setminus \bar H)$ and $G_{p,q}$.\qed

\section{Symplectomorphisms of generic coverings of the plane}
\label{Sect:Last}

An algebraic ramified covering of $\mathbb C \mathbb P^2$ is a
finite holomorphic map $f:X\to\mathbb C \mathbb P^2$ of a normal
projective irreducible complex surface $X$ to the projective
plane. The {\it ramification divisor\/} $R\subset X$ is the
divisor of the jacobian of $f$ (the multiplicity of $R$ is the
local degree of $f$ minus $1$). The {\it branch curve\/} $\bar
H\subset\mathbb C \mathbb P^2$ is the image of the support of $R$
or, in other words, the set of points over which $f$ is not
locally invertible.

The fundamental group $\pi_1=\pi_1(\mathbb C \mathbb P^2\setminus
\bar H,p)$ of the complement of $\bar H$, where $p\in \mathbb C
\mathbb P^2\setminus \bar H$, acts on the fibre $f^{-1}(p)$. Thus,
a homomorphism ({\it monodromy  of degree $N=\deg f$})
$\overline{\mu}=\overline{\mu} (f) : \pi_1 \to \mathfrak S_N$ from
$\pi_1$ to the symmetric group $\mathfrak S_N$ is well-defined.
The monodromy $\overline{\mu}$ is determined by $f$ uniquely up to
inner automorphism of the symmetric group. Conversely, by Grauert
- Remmert theorem (\cite{Gr-R}), a homomorphism $\overline{\mu} :
\pi_1 \to \mathfrak S_N$ the image $Im\, \overline{\mu}$ of which
acts transitively on a set consisting of $N$ elements is the
monodromy of some finite morphism $f:X\to \mathbb C \mathbb P^2$
of $\deg f=N$.

To describe the fundamental group of the complement of a curve
$\bar H$ of $\deg \bar H=m$, we fix a point $p\in \mathbb C
\mathbb P^2\setminus \bar H $, choose a point $x\in \bar
H\setminus Sing\, \bar H$, and consider a projective line $\Pi
=\mathbb C \mathbb P^1\subset \mathbb C \mathbb P^2$ meeting $\bar
H$ transversally at $x$. Let $\Gamma \subset \Pi$ be a circle of
small radius with center at $x$.  The complex orientation on
$\mathbb C \mathbb P^2$ defines an orientation on $\Gamma $. Let
$\gamma $ be a loop consisting of a path $L$ in $\mathbb C \mathbb
P^2\setminus \bar H$, joining $p$ with a point $q\in \Gamma$, the
loop $\Gamma $ (with positive direction) starting and ending at
$q$, and a return to $p$ along $L$ in the opposite direction. Such
loops $\gamma$ (and the elements in $\pi _1$ corresponding them)
are called {\it geometric generators} of the fundamental group
$\pi _1=\pi _1(\mathbb C \mathbb P^2\setminus \bar H,\, p)$. It is
well-known that $\pi _1$ is generated by the geometric generators,
and if $\bar H$ is an irreducible curve then any two geometric
generators are conjugated in $\pi _1$.

For each singular point $s_i$ of multiplicity $k_i$ of the curve
$\bar H$, let us choose a small neighborhood $U_i\subset \mathbb C
\mathbb P^2$ such that $\bar H \cap U_i$ is given (in local
analytic coordinates in $U_i$) by Weierstrass equation
$$w^{k_i}+\sum_{j=0}^{k_i-1}c_j(z)w^j=0$$
(in particular,  $w^2=z^3$ if $s_i$ is a cusp and $w^2=z^2$ if
$s_i$ is a node). Let $p_i$ be an arbitrary point in $U_i\setminus
\bar H$. Recall that if $s_i$ is a cusp then $\pi _1(U_i\setminus
\bar H,p_i)$ is isomorphic to the braid group $Br_3$ on three
strings and generated by two geometric generators (say, $a_1$ and
$a_2$) satisfying the relation $a_1a_2a_1=a_2a_1a_2$, and if $s_i$
is a node then $\pi _1(U_i\setminus \bar H,p_i)$ is isomorphic to
$\mathbb Z\oplus \mathbb Z$ and generated by two commuting
geometric generators.

Choose smooth paths $L_i$ in $\mathbb C \mathbb P^2\setminus \bar
H$, connecting the points $p_i$ with $p$. This choice defines
homomorphisms $\psi _i:\pi _1(U_i\setminus \bar H,p_i)\to \pi _1$.
We call  $\psi _i(\pi _1(U_i\setminus \bar H,p_i))=G_i$  {\it the
local fundamental group} of the singular point $s_i$. The local
fundamental groups $G_i$ depend on the choice of the paths $L_i$
and therefore are defined uniquely up to conjugation in $\pi_1$.
If $bmf(\bar H)=s=u_1\cdot ...\cdot u_n\in S_{Br_m}$ and the
factor $u_i$ corresponds to the singular point $s_i$ of $\bar H$,
then, in view of Zariski -- van Kampen Theorem, the group $G_i$
coincides (up to conjugation) with the image of $G_{u_i,loc}$
under the natural epimorphism $r:G(B_s)\to  G(B_s)/ \{
\varphi_{W_s}(x_1...x_m)=\mathbf 1 \}$. The monodromy
$\overline{\mu} :  \pi_1 \to \mathfrak S_N$ of a generic covering
$f:X\to \mathbb C\mathbb P^2$ branched along $\bar H$ defined a
homomorphism $\mu:G(s)\to \mathfrak S_N$, where $\mu
=\overline{\mu}\circ r$ and $s=bmf(\bar H)$, which we also call
the {\it monodromy} of $f$.

A finite morphism $f:X\to \mathbb C \mathbb P^2$, $\deg f=N$,
ramified over a cuspidal curve $\bar H$ is called {\it a generic
covering of the plane with discriminant curve} $\bar H$ if the
monodromy $\mu $ of $f$ is a generic epimorphism, i.e., it
satisfies the following conditions:
\begin{itemize}
\item[(i)] $\mu$ is an epimorphism,
\item[(ii)] the
image $\mu (\gamma)$ of each geometric generator $\gamma$ is a
transposition in $\mathfrak S_N$,
%and
\item[(iii)] the images $\mu (G_i)$ of
the local fundamental groups $G_i$ corresponding to the singular
points of $\bar H$ are subgroups of $\mathfrak S_N$ of orders
greater than 2.
\end{itemize}
The last condition is equivalent to the condition
that $X$ is a non-singular algebraic surface (see, for example,
\cite{Ku}). Note that if $\bar H$ is the discriminant curve of
some generic covering, then its degree is an even number.  For
a generic covering $f$, the total transform $f^*(\bar H)$ splits in
the sum $f^*(\bar H)=2R+C$, where $R\subset X$ is the ramification
divisor of $f$ having ramification multiplicity equal to two and
$C$ is a reduced curve in which $f$ is not ramified. The
restriction $f_{\mid R}$ to $R$ is a morphism of degree 1.

In the symplectic case, by an observation of Gromov, the
discussion of generic coverings of the plane given in terms of
monodromies of generic coverings can be generalized to the
symplectic situation as well (see Lemma 1 in \cite{Go} and
Proposition 10 in \cite{Au} for more details). Firstly, given a
cuspidal symplectic surface $\bar H$ (possibly with negative
nodes) in $\mathbb C \mathbb P^2$ and a generic epimorphism
(monodromy) $\overline{\mu}$ of the fundamental group of its
complement onto ${\mathfrak S}_N$, one can construct a smooth
ramified covering $f:X\to\mathbb C \mathbb P^2$ of degree $N$
branched over $\bar H$. The preimage of the symplectic form
vanishes on the normal bundle to the ramification surface
$R\subset X$ and at the preimages of cusps in $R$. However, it is
possible to add an exact $2$-form $\eta\in\Lambda^2(X)$ so that
$f^*\omega_{\mathbb C \mathbb P^2}+\eta$ is symplectic. A subtle
point is that, although there is no canonical way to lift the
symplectic {\it form\/} to $X$, the resulting symplectic {\it
structure\/} is defined canonically by a Moser-type argument.

We state the result in the form of a lemma.

\begin{lem}
\label{SympCov} Suppose that $\bar H\subset\mathbb C \mathbb P^2$
is a cuspidal symplectic surface (possibly with negative nodes)
and $\overline{\mu} :\pi_1(\mathbb C \mathbb P^2\setminus \bar
H)\to {\mathfrak S}_n$ is a generic monodromy. Then there exists a
unique symplectic four-manifold $X$ with a smooth map
$f:X\to\mathbb C \mathbb P^2$ branched over $\bar H$ as described
above such that the symplectic structure on $X$ is the pull-back
of the standard symplectic structure on $\mathbb C \mathbb P^2$.
\end{lem}

Clearly, given a symplectic isotopy $\bar H_t,\,t\in[0,1]$, of
cuspidal symplectic surfaces such that there exists an epimorphism
$\overline{\mu} :\pi_1(\mathbb C \mathbb P^2\setminus \bar H_0)\to
{\mathfrak S}_N$, we obtain a family of symplectomorphic
symplectic 4-manifolds $X_t$. (Recall that, by the definition of
isotopy, the topology of $\mathbb C \mathbb P^2\setminus \bar H_t$
remains unchangeable.) Applying similar arguments as in
\cite{Kh-Ku}, one can show that any cuspidal symplectic surface
with negative nodes in $\mathbb C\mathbb P^2$ is symplectically
isotopic to a cuspidal symplectic Hurwitz curve with negative
nodes. So in symplectic case to investigate generic covering of
$\mathbb C\mathbb P^2$ branched along a cuspidal symplectic
surface $\bar H$, we can restrict ourselves to the case when $\bar
H$ is a Hurwitz curve.

Consider two generic coverings $f^{\prime}:X^{\prime}\to \mathbb
C\mathbb P^2$ and $f^{\prime\prime}:X^{\prime\prime}\to \mathbb
C\mathbb P^2$ of degree $N$ branched, respectively, along cuspidal
Hurwitz curves (possibly with negative nodes) $\bar H^{\prime}$
and $\bar H^{\prime\prime}$. Let $s^{\prime}=bmf (\bar
H^{\prime})$ and $s^{\prime\prime}=bmf (\bar H^{\prime\prime})$ be
braid monodromy factorizations of $\bar H^{\prime}$ and $\bar
H^{\prime\prime}$, and $f^{\prime}$ (respectively,
$f^{\prime\prime}$) defined by the monodromy
$\mu^{\prime}:G(s^{\prime})\to \mathfrak S_N$ (respectively,
$\mu^{\prime\prime}$) We say that $\bar H^{\prime}$ and $\bar
H^{\prime\prime}$ are {\it weakly $\mu$-equivalent} if the pairs
$(s^{\prime},\mu_{s^{\prime}})$ and
$(s^{\prime\prime},\mu_{s^{\prime\prime}})$ are weakly
$\mu$-equivalent.

\begin{lem} \label{le}
Let two Hurwitz curves $\bar H^{\prime}$ and $\bar
H^{\prime\prime}$ be the branch curves of two generic coverings
$f^{\prime}$ and $f^{\prime\prime}$ of the plane. If $\bar
H^{\prime}$ and $\bar H^{\prime\prime}$ are weakly
$\mu$-equivalent, then there is a symplectomorphism
$h:X^{\prime}\to X^{\prime\prime}$ such that the coverings
$f^{\prime}$ and $f^{\prime\prime}\circ h$ are regular homotopic.
\end{lem}
\proof It follows from the definition of weak $\mu$-equivalence
Claim \ref{hurweak},
and Theorem 3.3 in \cite{Kh-Ku}. \qed

In algebraic case, if $(X,L)$ is a polarized projective surface,
that is, $L$ is an ample line bundle on $X$, then, for
sufficiently large $k$, the sections of $L^{\otimes k}$ define an
embedding of $X$ into some projective space  $\mathbb C \mathbb
P^r$. The restriction of a generic projection $\mathbb C \mathbb
P^r\to\mathbb C \mathbb P^2$ gives rise to a generic ramified
covering $f_k:X\to\mathbb C \mathbb P^2$ with branch curve $\bar
H_k$. Observe that the degree of this covering equals $k^2\deg
L=k^2c_1(L)^2$. The space of generic projections is path-connected
and therefore we can give the following definition.

\begin{df}
\label{SurfInv} The $k$-th braid monodromy invariant $\mu_k(X,L)$
of the polarized complex surface $(X,L)$ is the type of the braid
monodromy factorization of $\bar H_k$.
\end{df}

The polarization canonically defines a symplectic structure
$\omega_L$ on $X$. To see this, pick a positive $(1,1)$-form
representing the first Chern class $c_1(L)\in H^2(X,\mathbb Z)$.
Every two such forms can be joined by a linear homotopy, hence
they are diffeomorphic by Moser's theorem.

\begin{thm} \label{t}
Let $(X,L)$ and $(X',L')$ be polarized projective surfaces.
Suppose that $\mu_k(X,L)=\mu_k(X',L')$ for some $k$ such that
$k^2c_1(L)^2\ge 12$. Then $(X,\omega_L)$ and $(X',\omega_{L'})$
are symplectomorphic.
\end{thm}

\proof By assumption, the branch curves $\bar H$ and $\bar H'$ of
generic projections onto $\mathbb C \mathbb P^2$, defined by
sections of $L^{\otimes k}$ and $(L')^{\otimes k}$, respectively,
have the same type of braid monodromy factorizations. Thus, by
Theorem 3.3 in \cite{Kh-Ku}, $\bar H$ and $\bar H'$ are
symplectically isotopic as cuspidal Hurwitz curves. By virtue of
Theorem 3 in \cite{Ku2}, there exists a unique generic monodromy
$\mu :\pi_1(\mathbb C \mathbb P^2\setminus \bar H)\to {\mathfrak
S}_N$ for $N=k^2c_1(L_1)^2$. Therefore, by Lemma \ref{SympCov}, we
have a unique family of symplectic ramified coverings of $\mathbb
C\mathbb P^2$ connecting $X$ and $X'$. It follows from Moser's
theorem that $X$ and $X'$ are symplectomorphic with the pulled
back symplectic structures. However, these pull-backs are
proportional to the symplectic structures defined by the
polarizations with the same coefficient $k$, which completes the
proof.\qed

\smallskip

Let $(X,\omega)$ be a compact symplectic 4-manifold with
symplectic form $\omega$ whose class $[\omega]\in H^2(X,\mathbb
Z)$. Fix an $\omega$-compatible almost complex structure $J$ and
corresponding Riemannian metric $g$. Let $L$ be a line bundle on
$X$ whose first Chern class is $[\omega]$. By \cite{Au}, for
$k>>0$, the line bundle $L^{\otimes k}$ admits many approximately
holomorphic sections so that one can choose three of them which
give an approximately holomorphic generic covering $f_k:X\to
\mathbb C\mathbb P^2$ of degree $N_k=k^2 \omega^2$ branched alone
a cuspidal Hurwitz curve $\bar H$ (possibly with negative nodes).
Denote, as above, by $bmf(\bar H_k)=s_k$ a braid monodromy
factorization of the Hurwitz curve $\bar H_k$ and $\mu_k
:G(B_{s_k})\to \mathfrak S_{N_k}$ the monodromy of $f_k$.
\begin{df}
\label{SymInv} The $k$-th braid monodromy invariant
$\overline{\mu}_k(X,\omega)$ of a symplectic 4-manifold
$(X,\omega)$ is the pair $(bmf(\bar H_k), \mu_k)$.
\end{df}

\begin{thm} {\rm ( \cite{Au-Ka})}
If two symplectic 4-manifolds $(X^{\prime},\omega^{\prime})$ and
$(X^{\prime\prime},\omega^{\prime\prime})$ are symplectomorphic,
then, for $k$ sufficiently large, the braid monodromy invariants
$\mu_k(X^{\prime},\omega^{\prime})$ and
$\mu_k(X^{\prime\prime},\omega^{\prime\prime})$  are weakly
$\mu$-equi\-va\-lent.
\end{thm}

Weak equivalence appears in this statement, since the construction
of an isotopy of the coverings depends continuously on a choice of
tamed almost complex structures. Therefore, for symplectomorphic
4-manifolds, we obtain a "regular homotopy" of branch curves,
which leads to the weak $\mu$-equivalence of their braid monodromy
factorizations. The following statement shows that in general case
we can not avoid the appearance of negative nodes.

\begin{thm} \label{last} For any $N\geq 4$, there
are two isotopic generic coverings
$f_1: X\to \mathbb C\mathbb P^2$ and $f_2: X\to \mathbb C\mathbb
P^2$ of degree $N$ branched along cuspidal Hurwitz curves $\bar
H_1$ and $\bar H_2$ without negative nodes such that
\begin{itemize}
\item[(i)] $\bar H_1$ and $\bar H_2$ have different braid
monodromy factorization types;
\item[(ii)] The pairs $(bmf(\bar H_1), \mu_1)$ and $(bmf(\bar H_2),
\mu_2)$ are weakly $\mu$-equi\-va\-lent.
\end{itemize}
\end{thm}
\proof Let $\bar H_N$ be the branch curve of a generic linear
projection $\mbox{pr}:X_N\to \mathbb C\mathbb P^2$ of a smooth
projective surface $X_N\subset \mathbb C\mathbb P^3$ of $\deg
X_N=N$. Denote by $m=N(N-1)$ the degree of $\bar H_N$ and
$s=s_N=bmf(\bar H_N)$ its braid monodromy factorization. The
projection $\mbox{pr}$ defines a generic epimorphism (monodromy)
$\mu=\mu_N :G(B_s)\to \mathfrak S_N$. By \cite{Moi1}, the
$C$-group
\begin{equation}
\begin{array}{l}
G(s)=G(W(B_s,\mathbb F_m))\simeq \\
<x_1, \dots , x_m \, \mid \, x_i^{-1}b(x_i)=\mathbf 1, \, \,
i=1,\dots, m,\, \, b\in B_s\,
>
\end{array}
\end{equation}
is $C$-isomorphic to the braid group $Br_N$ and $\mu$ coincides
with canonical epimorphism $\sigma :Br_N\to \mathfrak S_N$ with
the group of pure braids as its kernel. Without loss of
generality, we can assume that $\varphi _s(x_1)=a_1\in Br_N$.
Since $\mu$ is an epimorphism, there is an element $y$ conjugated
in $\mathbb F_m$ to some generator $x_i$ such that
$\varphi_s(y)=a_2^{-2}a_3a_2^2$. Note that $[\varphi
_s(x_1),\varphi _s(y)]\neq \mathbf 1$ in $Br_N$, but $\mu(\varphi
_s(x_1))=(1,2)$ and $\mu(\varphi _s(y))=(3,4)$ are two different
commuting transpositions in $\mathfrak S_N$. Then, by Theorem
\ref{alg}, there are $M\in \mathbb N$ and two pairs $(\overline
s,\mu_{\overline s})$ and $(\widetilde s,\mu_{\widetilde s})$,
$\overline s,\, \widetilde s\in \mathcal A_2\subset S_{Br_M}$,
such that
\begin{itemize}
\item[(1)] $\alpha (\overline s)=\alpha (\widetilde s)=\Delta^2_M$,
\item[(2)] $G_{\overline s}$ is $C$-isomorphic to $G(s)\simeq Br_N$ and
$G_{\widetilde s}$ is $C$-isomorphic to the $C$-group $Br_N/\{
[a_1,a_2^{-2}a_3a_2^2]=\mathbf 1\}$,
\item[(3)]  the types of $\overline s$ and $\widetilde s$ are different,
\item[(4)] $\mu_{\overline s}$ and $\mu_{\widetilde s}$
are generic epimorphisms onto $\mathfrak S_N$,
\item[(5)] the pairs $(\overline s,\mu_{\overline s})$ and
$(\widetilde s,\mu_{\widetilde s})$ are weakly $\mu$-equivalent
\item[(6)] $\overline s,\, \widetilde s\in \mathcal A^0_2\subset
S_{Br_M}$.
\end{itemize}
Now by Theorem \ref{moisha}, there are two cuspidal Hurwitz curves
$\bar H_1$ and $\bar H_2$ in $\mathbb C\mathbb P^2$ of degree $M$
with $bmf(\bar H_1)=\overline s$ and $bmf(\bar H_2)=\widetilde s$.
The generic epimorphisms $\mu_{\overline s}$ and $\mu_{\widetilde
s}$ define two generic coverings $f^{\prime}:X^{\prime}\to \mathbb
C\mathbb P^2$ and $f^{\prime\prime}:X^{\prime\prime}\to \mathbb
C\mathbb P^2$ branched, respectively, along $\bar H_1$ and $\bar
H_2$. By Lemma \ref{le}, there is a symplectomorphism
$h:X^{\prime}\to X^{\prime\prime}$ such that the coverings
$f_1=f^{\prime}:X=X^{\prime}\to \mathbb C\mathbb P^2$ and
$f_2=f^{\prime\prime}\circ h:X\to \mathbb C\mathbb P^2$ are
isotopic, but their branch Hurwitz curves have different braid
monodromy factorization types.\qed

\ifx\undefined\bysame
\newcommand{\bysame}{\leavevmode\hbox to3em{\hrulefill}\,}
\fi

\end{document}